\documentclass[pdflatex,sn-mathphys-num]{sn-jnl}
\usepackage{graphicx}
\usepackage{amsmath,amssymb,amsfonts,amsthm}
\usepackage{bm}
\usepackage{booktabs,multirow,makecell}
\usepackage{enumitem,placeins}
\usepackage{etoolbox}
\makeatletter
\patchcmd{\@maketitle}{\removelastskip\vskip36pt\vskip0pt}{%
  \par\medskip{\keywordfont\noindent\textbf{Acknowledgments}\quad
  Liying Zhang is supported by the National Natural Science Foundation of China
  under \mbox{Grant No.~12671464}.\par}%
  \removelastskip\vskip36pt\vskip0pt
}{}{\PackageError{MGRIT-BIT}{Title-page acknowledgment could not be inserted}{Check the official template version.}}
\makeatother

\newcommand{\R}{\mathbb{R}}
\newcommand{\N}{\mathbb{N}}
\newcommand{\E}{\mathbb{E}}

\theoremstyle{thmstyleone}
\newtheorem{theorem}{Theorem}[section]
\newtheorem{lemma}{Lemma}

\newtheorem{cor}{Corollary}
\theoremstyle{thmstylethree}

\newtheorem{remark}{Remark}

\newtheorem{assumption}{Assumption}
\numberwithin{equation}{section}
\allowdisplaybreaks[1]
\hypersetup{
  pdftitle={Mean-Square Convergence Analysis of a Three-Level MGRIT Method for Linear Stochastic Evolution Equations with Additive Noise},
  pdfauthor={Liying Zhang and Zhiheng Xu},
  pdfkeywords={three-level MGRIT, F-relaxation, stochastic evolution equations, mean-square error estimates, convergence analysis}
}
\begin{document}
\title[Mean-square convergence of three-level MGRIT]{Mean-Square Convergence Analysis of a Three-Level MGRIT Method for Linear Stochastic Evolution Equations with Additive Noise}
% Corresponding author designation requires author confirmation.
\author[1]{\fnm{Liying} \sur{Zhang}}\email{lyzhang@lsec.cc.ac.cn}
\author[1]{\fnm{Zhiheng} \sur{Xu}}\email{2528274616@qq.com}
\affil[1]{\orgdiv{School of Mathematical Science}, \orgname{China University of Mining and Technology (Beijing)}, \orgaddress{\city{Beijing}, \postcode{100083}, \country{China}}}
\abstract{{\unboldmath We investigate the mean-square convergence of a three-level multigrid reduction-in-time (MGRIT) method based on F-relaxation for linear stochastic evolution equations driven by additive Gaussian noise. The method combines exponential Euler time stepping on the finest grid with implicit Euler discretizations on the two coarser grids. The inexact solution of the intermediate-level full approximation scheme problem gives rise to coupled error recursions involving two interlevel propagator discrepancies. Under a common spectral decomposition of the drift and noise covariance operators, we establish a geometric contraction estimate in a modal mean-square norm, with a contraction factor uniform in the temporal grids and both coarsening factors. We further derive fixed-iteration temporal error bounds that are uniform in the first coarsening factor and explicitly quantify the dependence on the second.
	For a fixed iteration count  $k$ and a fixed second coarsening factor, we prove an $O(\Delta T_1^r)$ error bound relative to the finest-grid reference solution for $0<r\le k+1$ provided that $\mathcal R_{r,B}<\infty$.	
		 We also prove that the error-propagation operator is nilpotent on each fixed temporal grid, ensuring finite termination of the iteration. Numerical experiments illustrate the theoretical convergence behavior and compare temporal critical-path work under a unit-cost model.
}}

\maketitle

\section{Introduction}
The numerical solution of evolution equations commonly relies on sequential time integration, in which each time step depends on information from preceding steps. This temporal dependence can limit computational scalability when spatial parallelism is insufficient or has already been exploited. Parallel-in-time (PinT) methods seek to alleviate this limitation by enabling concurrent computations over different time intervals. For stochastic evolution equations, however, temporal parallelization must also respect the coupling of stochastic increments across time grids. Moreover, the attainable mean-square accuracy depends on the spatial regularity of the driving noise, making the analysis of temporal multilevel iterations more delicate than in the deterministic setting.

The parareal algorithm of Lions, Maday, and Turinici~\cite{LionsMadayTurinici2001} introduced an iterative coupling of fine and coarse time propagators; its convergence properties were subsequently studied by Gander and Vandewalle~\cite{GanderVandewalle2007}. Multigrid reduction in time (MGRIT) embeds temporal parallelism in a multigrid framework through relaxation and coarse-grid correction~\cite{Falgout2014,Falgout2017}. At coarse-grid points, the standard two-level MGRIT method with F-relaxation is equivalent to parareal. Deterministic convergence analyses of two-level and multilevel MGRIT have been developed in~\cite{DobrevEtAl2017,HessenthalerEtAl2020}. Although these results provide a foundation for understanding multilevel temporal corrections, their application to stochastic evolution equations requires additional control of noise-driven errors and their spectral regularity.

For stochastic partial differential equations, the regularity of the stochastic forcing can constrain the strong convergence rates of time discretizations; see, for example, the implicit Euler analysis in~\cite{Printems2001}. Br\'ehier and Wang~\cite{BrehierWang2020} analyzed parareal methods for semilinear parabolic stochastic partial differential equations and, in the linear additive-noise setting, obtained error bounds relative to a fine-grid reference solution for exponential Euler and linear implicit Euler propagators. Their results relate the attainable temporal convergence orders to both the iteration count and noise regularity. Related developments include randomized parareal methods for deterministic ordinary differential equations~\cite{PentlandEtAl2023,PentlandTamborrinoSullivan2023}, micro--macro parareal methods for McKean--Vlasov equations~\cite{BossuytVandewalleSamaey2025}, and parareal methods for stochastic Maxwell equations~\cite{ZhangZhang2024}. These studies provide important
theoretical and methodological
foundations for stochastic
parallel-in-time integration.
Building on this foundation,
the three-level MGRIT framework
introduces an additional
coarse-grid correction and its
associated interlevel coupling.
Analyzing this coupling, together
with the dependence on the second
coarsening factor, calls for a
corresponding extension of the
two-level error analysis.

The main analytical difficulty is the coupling induced by the inexact solution of the intermediate-level full approximation scheme (FAS) problem. In the two-level parareal framework,
the error analysis involves the
discrepancy between the fine and
coarse propagators. In the three-level
setting considered here, an additional
discrepancy arises between repeated
intermediate-level propagation and
propagation on the coarsest level.
The resulting coupled error structure
requires the simultaneous treatment
of both interlevel discrepancies.
To account for this additional
correction, we introduce an
auxiliary error associated with
the exact intermediate-level FAS
solve and analyze its coupling
with the error induced by the
coarsest-level correction.
 A further difficulty arises from the unbounded spectrum of the drift operator. Uniform operator bounds alone do not capture the frequency-dependent decay needed to obtain positive powers of the temporal step size after summing over infinitely many stochastic modes. Accordingly, the analysis must combine uniform bounds for normalized interlevel propagator discrepancies with estimates that retain the modal time scales and exploit spectral summability of the noise.

We study a three-level MGRIT method with F-relaxation (F-MGRIT) for the linear stochastic evolution equation
\begin{equation*}
	\mathrm{d}u(t)=(A+B)u(t)\,\mathrm{d}t+Q^{1/2}\,\mathrm{d}W(t),
	\qquad u(0)=0,
\end{equation*}
posed on a real separable Hilbert space $H$. Here $A$ is a negative definite self-adjoint operator, $B$ is a bounded self-adjoint perturbation, $Q$ is a nonnegative self-adjoint noise covariance operator, and $W$ is a cylindrical Wiener process. We assume that $A$, $B$, and $Q$ have a common orthonormal eigenbasis and that $-(A+B)$ has a strictly positive spectral lower bound; the precise assumptions are stated in Section~\ref{sec:setting}. The finest temporal grid (Level 0) uses exponential Euler time stepping, while the two coarser grids (Levels 1 and 2) use linear implicit Euler discretizations. The grids are nested with coarsening factors $J_1$ and $J_2$, and coarse-grid Brownian increments are formed by summing the corresponding fine-grid increments. All iteration errors are measured relative to the Level-0 reference trajectory restricted to Level-1 time points. By resolving the coupled errors mode by mode, we obtain both a coarsening-uniform geometric bound and frequency-sensitive estimates for temporal refinement.

The main contributions are as follows.
\begin{enumerate}
	\item  We derive Hilbert-space and modal error recursions that separate the contribution of the Level-0--Level-1 propagator discrepancy from the additional correction error associated with the Level-1--Level-2 discrepancy. This formulation explicitly accounts for the inexact intermediate-level FAS solve and reduces to the corresponding two-level parareal recursion when that solve is exact.
	
	\item We establish a geometric error bound with a contraction factor $q_\ast<1$ that is independent of the temporal grids and both coarsening factors. More precisely, if $X_{1,n}^{(k)}$ is the $k$th Level-1 iterate and $X_{1,n}^{\mathrm{ref}}$ is the Level-0 reference solution restricted to the Level-1 grid, then for $k\in\mathbb N_0$ and $0<r\le k+1$ satisfying $\mathcal R_{r,B}<\infty$, we prove
	\begin{equation*}
		\max_{0\le n\le N_1}
		\bigl\|X_{1,n}^{(k)}-X_{1,n}^{\mathrm{ref}}\bigr\|_{L^2(\Omega;H)}
		\le C_{\mathrm{init}}q_\ast^k J_2^{(r-1)_+}
		\mathcal R_{r,B}\,\Delta T_1^r,
	\end{equation*}
	where $\mathcal R_{r,B}:=\bigl(\sum_{p\ge1}\eta_p^{2r-1}\gamma_p\bigr)^{1/2}$, $(r-1)_+:=\max\{r-1,0\}$, and $\eta_p$ and $\gamma_p$ are the eigenvalues of $-(A+B)$ and $Q$, respectively. The constant $C_{\mathrm{init}}$ is independent of the iteration count and temporal grid parameters. Thus, for fixed $J_2$, the iteration error is $O(\Delta T_1^r)$; for $0<r\le1$, the bound is uniform in both coarsening factors. We also identify the admissible orders for power-law spectra and recover the regularity-dependent temporal orders of the corresponding two-level analysis in the one-dimensional Dirichlet Laplacian setting~\cite{BrehierWang2020}.
	
	\item We show that the iteration error-propagation operator is strictly lower block triangular and therefore nilpotent. Consequently, on every fixed temporal grid, the Level-0 reference values at all Level-1 points are recovered after at most $N_1$ outer iterations.
\end{enumerate}

Numerical experiments illustrate the predicted convergence behavior and compare temporal critical-path work for different PinT configurations under a unit-cost model. The geometric iteration bound, finite termination, and fixed-iteration temporal error estimates describe distinct aspects of convergence, all relative to the Level-0 reference trajectory; measured parallel runtimes and speedups are not addressed here.

The remainder of the paper is organized as follows. Section~\ref{sec:setting} introduces the stochastic evolution equation and spectral assumptions. Section~\ref{sec:discretization} describes the temporal discretizations and three-level F-MGRIT method. Section~\ref{sec:modal} establishes the coupled error recursions, geometric bounds, and fixed-grid finite termination. Section~\ref{sec:refinement} develops fixed-iteration mean-square error estimates under temporal refinement. Section~\ref{sec:numerical} presents numerical experiments, and Section~\ref{sec:conclusion} concludes the paper.

\section{Problem setting and preliminaries}
\label{sec:setting}
Let $H$ be a real separable Hilbert space with inner product
$(\cdot,\cdot)_H$ and norm $\|\cdot\|_H$. Let
$(\Omega,\mathcal F,\mathbb P,\{\mathcal F_t\}_{t\ge0})$
be a complete filtered probability space satisfying the
usual conditions. Fix $T>0$, and set
$\N:=\{1,2,\ldots\}$ and $\N_0:=\N\cup\{0\}$.
We write $a\wedge b:=\min\{a,b\}$ and use
$v_p\asymp w_p$ to indicate that
$cw_p\le v_p\le Cw_p$ for all sufficiently large $p$,
with constants $c,C>0$ independent of $p$.

We denote by $\mathcal L(H)$ and $\mathcal L_2(H)$
the spaces of bounded linear and Hilbert--Schmidt
operators on $H$, respectively. The identity operator
is denoted by $I$, and $\sigma(\cdot)$ denotes the
spectrum. For $Y\in L^2(\Omega;H)$, we set $\|Y\|_{L^2(\Omega;H)}
:=\bigl(\E\|Y\|_H^2\bigr)^{1/2},$
and write $L^2(\Omega,\mathcal F_t;H)$ for the
subspace of $\mathcal F_t$-measurable random
variables. Throughout the paper, $C>0$ denotes
a generic constant that may change from line to
line. Unless otherwise stated, it is independent
of the temporal discretization parameters and
modal indices.

Set $A_B:=A+B$, with $D(A_B):=D(A)$.
We consider the linear stochastic evolution equation
\begin{equation}
	\mathrm du(t)
	=A_Bu(t)\,\mathrm dt
	+Q^{1/2}\,\mathrm dW(t),
	\qquad u(0)=0,
	\label{eq:model}
\end{equation}
where $W$ is a cylindrical Wiener process on $H$.
The operator $Q$ specifies the spectral covariance
structure of the additive Gaussian forcing.
The stochastic term in~\eqref{eq:model} is
understood through the stochastic convolution
constructed below; in particular, $Q$ is not
required to be trace class.

To facilitate the modal analysis of the temporal
propagators, we impose the following common
spectral structure.
\begin{assumption}
	\label{ass:operators}
	The operator $A:D(A)\subset H\to H$ is
	self-adjoint and negative definite,
	$B\in\mathcal L(H)$ is self-adjoint, and
	$Q:D(Q)\subset H\to H$ is nonnegative and
	self-adjoint. There exists a complete
	orthonormal basis
	$\{e_p\}_{p\ge1}\subset D(A)\cap D(Q)$
	consisting of common eigenvectors of
	$A$, $B$, and $Q$. For every $p\ge1$, $	Ae_p=-\lambda_pe_p,
	Be_p=\mu_pe_p,$ and
	$Qe_p=\gamma_pe_p,$
	where $\lambda_p,\mu_p,\gamma_p\in\mathbb R$
	and $\gamma_p\ge0$. Moreover, there exists
	$\eta_0>0$ such that
	$
	\eta_p:=\lambda_p-\mu_p\ge\eta_0, p\ge1.
$
\end{assumption}
Since $B$ is bounded and self-adjoint,
the bounded perturbation theorem implies
that $A_B$ is self-adjoint on $D(A)$;
see~\cite{Pazy1983}. Define $L_B:=-A_B$.
Assumption~\ref{ass:operators} gives $
A_Be_p=-\eta_pe_p,$ and 
$L_Be_p=\eta_pe_p, p\ge1.$
Consequently, $
L_B\ge\eta_0I,
\sigma(A_B)\subset(-\infty,-\eta_0].
$
In particular, $A_B$ generates an analytic
$C_0$-semigroup of contractions
$\{\mathrm e^{tA_B}\}_{t\ge0}$ satisfying
$
\|\mathrm e^{tA_B}\|_{\mathcal L(H)}
\le \mathrm e^{-\eta_0t},
 t\ge0.
$
Furthermore, for every $\tau>0$,
\[
\|(I-\tau A_B)^{-1}\|_{\mathcal L(H)}
\le(1+\eta_0\tau)^{-1}.
\]
These bounds provide the stability properties
of the exponential and implicit Euler
propagators used below.

We next introduce the fractional powers of
$L_B$ and the associated spectral regularity
conditions. For $s\ge0$, the spectral calculus
defines
\[
D(L_B^s)
:=
\left\{
v\in H:
\sum_{p\ge1}\eta_p^{2s}
|(v,e_p)_H|^2<\infty
\right\},
\]
and $
L_B^sv
=
\sum_{p\ge1}
\eta_p^s(v,e_p)_He_p,
 v\in D(L_B^s).
$
For $s<0$, the same spectral expansion defines
a bounded operator on $H$, with
$
\|L_B^s\|_{\mathcal L(H)}
\le\eta_0^s.$
Similarly,
\[
D(Q^{1/2})
=
\left\{
v\in H:
\sum_{p\ge1}\gamma_p
|(v,e_p)_H|^2<\infty
\right\},
\]
and
$
Q^{1/2}v
=
\sum_{p\ge1}
\sqrt{\gamma_p}(v,e_p)_He_p,
 v\in D(Q^{1/2}).
$
For any $s\in\mathbb R$, the spectral
operator $L_B^sQ^{1/2}$, initially defined
on $\operatorname{span}\{e_p:p\ge1\}$,
admits a unique Hilbert--Schmidt extension
to $H$ whenever $
\sum_{p\ge1}\eta_p^{2s}\gamma_p<\infty.
$
In this case,
$
\|L_B^sQ^{1/2}\|_{\mathcal L_2(H)}^2
=
\sum_{p\ge1}\eta_p^{2s}\gamma_p.
$
The spectral calculus and fractional-power
construction are standard; see~\cite{Pazy1983}.

Let $\{\beta_p\}_{p\ge1}$ be independent
standard Brownian motions adapted to
$\{\mathcal F_t\}_{t\ge0}$, with increments
after time $s$ independent of $\mathcal F_s$.
The cylindrical Wiener process is represented
formally by $
W(t)=\sum_{p\ge1}\beta_p(t)e_p.
$
Accordingly, the formal noise process has
the expansion
$
Q^{1/2}W(t)
=
\sum_{p\ge1}
\sqrt{\gamma_p}\,\beta_p(t)e_p.
$
If $\operatorname{Tr}(Q)
=\sum_{p\ge1}\gamma_p<\infty$,
this series defines an $H$-valued
$Q$-Wiener process. We do not impose
this trace-class condition. In general,
the preceding series is only formal,
and the stochastic forcing is defined
through the semigroup-smoothed stochastic
integral.

For $r\ge0$, define
$\mathcal R_{r,B}
	:=
	\left(
	\sum_{p\ge1}
	\eta_p^{2r-1}\gamma_p
	\right)^{1/2}
	\in[0,\infty].$
Whenever $\mathcal R_{r,B}<\infty$,
the operator $L_B^{r-1/2}Q^{1/2}$
has a Hilbert--Schmidt extension and
$
\mathcal R_{r,B}
=
\|L_B^{r-1/2}Q^{1/2}\|_{\mathcal L_2(H)}.
$
The case $r=0$ provides the basic
summability condition required for
the stochastic convolution.

\begin{assumption}
	\label{ass:noise}
	The noise satisfies $
	\mathcal R_{0,B}^2
	=
	\sum_{p\ge1}
	\frac{\gamma_p}{\eta_p}
	<\infty.
$
\end{assumption}

We now verify the well-posedness of the
stochastic convolution. Under
Assumptions~\ref{ass:operators}--\ref{ass:noise},
the operator
$\mathrm e^{(t-s)A_B}Q^{1/2}$,
initially defined on the finite span
of the common eigenvectors, admits a
Hilbert--Schmidt extension for every
$0\le s<t\le T$. Moreover,
\[
\begin{aligned}
	\int_0^t
	\bigl\|
	\mathrm e^{(t-s)A_B}Q^{1/2}
	\bigr\|_{\mathcal L_2(H)}^2
	\,\mathrm ds=
	\sum_{p\ge1}\gamma_p
	\int_0^t
	\mathrm e^{-2\eta_p(t-s)}
	\,\mathrm ds\le
	\frac12\mathcal R_{0,B}^2
	<\infty.
\end{aligned}
\]
Thus the stochastic convolution $u(t)
:=
\int_0^t
\mathrm e^{(t-s)A_B}Q^{1/2}
\,\mathrm dW(s)$
is well defined as an $H$-valued
It\^o integral. Equivalently,
\[
u(t)
=
\sum_{p\ge1}
\sqrt{\gamma_p}
\left(
\int_0^t
\mathrm e^{-\eta_p(t-s)}
\,\mathrm d\beta_p(s)
\right)e_p,
\]
where the series converges in
$L^2(\Omega;H)$.
The stochastic convolution belongs to
$C([0,T];L^2(\Omega;H))$ and is the
unique adapted mild solution of
\eqref{eq:model} in this class;
see~\cite{DaPratoZabczyk2014,
	LordPowellShardlow2014}.
For each $p\ge1$, its modal component
$u_p(t):=(u(t),e_p)_H$ satisfies
\[
\mathrm du_p(t)
=
-\eta_pu_p(t)\,\mathrm dt
+\sqrt{\gamma_p}\,\mathrm d\beta_p(t),
\qquad u_p(0)=0.
\]

Finally, since $\eta_p\ge\eta_0>0$,
for every $r>0$, $\mathcal R_{0,B}
\le\eta_0^{-r}\mathcal R_{r,B}.$
Hence $\mathcal R_{r,B}<\infty$ implies
Assumption~\ref{ass:noise}. The basic
condition $\mathcal R_{0,B}<\infty$
ensures mean-square well-posedness,
whereas the stronger conditions
$\mathcal R_{r,B}<\infty$ quantify the
spectral regularity used in the
fixed-iteration temporal error analysis.
The common eigenbasis also allows the
discrete propagators and their errors
to be analyzed mode by mode.

\section{Three-level temporal discretization and F-MGRIT iteration}
\label{sec:discretization}

We discretize~\eqref{eq:model} on three nested temporal grids, with Brownian increments consistently coupled across all levels. An exponential Euler discretization on Level~0 defines the reference trajectory, whereas implicit Euler propagators on Levels~1 and~2 are used to construct the coarse-grid corrections. Our analysis concerns the iteration error relative to the Level-0 reference trajectory. Throughout Sections~\ref{sec:discretization} and~\ref{sec:modal}, Assumptions~\ref{ass:operators} and~\ref{ass:noise} are in force.

Let $J_1,J_2\in\N$ be the coarsening factors from Level~0 to Level~1 and from Level~1 to Level~2, respectively, with $J_1,J_2\ge2$. Given $N_2\in\N$, set $N_1:=J_2N_2$ and $N_0:=J_1N_1$, and define
$
\delta t_0:=\frac{T}{N_0},
\Delta T_1:=\frac{T}{N_1},
\Delta T_2:=\frac{T}{N_2}.
$
Thus $\Delta T_1=J_1\delta t_0$, $\Delta T_2=J_2\Delta T_1$, and $T=N_1\Delta T_1=N_2\Delta T_2$. Write $t_n:=n\Delta T_1$ for $0\le n\le N_1$. Within each interval $[t_n,t_{n+1}]$, define the Level-0 points $t_{n,j}:=t_n+j\delta t_0$ for $0\le n<N_1$ and $0\le j\le J_1$, and set $t_{N_1,0}:=T$. The Level-2 points are $t_{mJ_2}$, $0\le m\le N_2$. Relative to the Level-2 grid, the Level-1 points $t_{mJ_2}$ are called \emph{C-points}, while $t_{mJ_2+s}$, $0\le m<N_2$ and $1\le s<J_2$, are called \emph{F-points}.

For $0\le n<N_1$, $1\le j\le J_1$, and $p\ge1$, define the Level-0 Brownian increments by
$\Delta\beta_{n,j,p}^{(0)}:=\beta_p(t_{n,j})-\beta_p(t_{n,j-1}).$
The increments on the coarser levels are obtained by aggregation:
\[
\begin{aligned}
	\Delta\beta_{n,p}^{(1)}
	&:=\sum_{j=1}^{J_1}\Delta\beta_{n,j,p}^{(0)}
	=\beta_p(t_{n+1})-\beta_p(t_n),\qquad 0\le n<N_1,\\
	\Delta\beta_{m,p}^{(2)}
	&:=\sum_{q=0}^{J_2-1}\Delta\beta_{mJ_2+q,p}^{(1)}
	=\beta_p(t_{(m+1)J_2})-\beta_p(t_{mJ_2}),\qquad 0\le m<N_2.
\end{aligned}
\]
Level-0 increments belonging to different modes or disjoint temporal intervals are independent, and
$\E\Delta\beta_{n,j,p}^{(0)}=0$ and
$\E|\Delta\beta_{n,j,p}^{(0)}|^2=\delta t_0$.

The deterministic parts of the one-step propagators are $
E_0:=\mathrm e^{\delta t_0A_B},$ 
$S_1:=(I-\Delta T_1A_B)^{-1},$ and
$S_2:=(I-\Delta T_2A_B)^{-1}.
$
We use exponential Euler on Level~0 and implicit Euler on Levels~1 and~2. Define the corresponding discrete noise increments by $
	\Delta\widetilde W_{n,j}^{(0)}
	:=\sum_{p\ge1}\sqrt{\gamma_p}\,\mathrm e^{-\eta_p\delta t_0}
	\Delta\beta_{n,j,p}^{(0)}e_p,$
	$\Delta\widetilde W_n^{(1)}
	:=\sum_{p\ge1}\sqrt{\gamma_p}\,
	\frac{\Delta\beta_{n,p}^{(1)}}{1+\eta_p\Delta T_1}e_p,$ and
$	\Delta\widetilde W_m^{(2)}
	:=\sum_{p\ge1}\sqrt{\gamma_p}\,
	\frac{\Delta\beta_{m,p}^{(2)}}{1+\eta_p\Delta T_2}e_p.$
For $h,\eta>0$, $
h\mathrm e^{-2\eta h}\le\frac{1}{2\mathrm e\eta},
\frac{h}{(1+\eta h)^2}\le\frac{1}{4\eta}.
$
Consequently, independence of the modal Brownian increments, It\^o's isometry, and Assumption~\ref{ass:noise} show that all three series converge in $L^2(\Omega;H)$. The Level-0 exponential Euler step approximates the stochastic convolution by an endpoint-weighted Brownian increment; it is not the exact Ornstein--Uhlenbeck transition, whose noise term contains a stochastic integral with a time-dependent exponential kernel.

For $Y\in L^2(\Omega;H)$, define the affine one-step propagators
$\Phi_{0,n,j}(Y):=E_0Y+\Delta\widetilde W_{n,j}^{(0)},$ 
$G_{1,n}(Y):=S_1Y+\Delta\widetilde W_n^{(1)},$ and
$	G_{2,m}(Y):=S_2Y+\Delta\widetilde W_m^{(2)}.$
Over one Level-1 interval, let $
F_{0,n}:=\Phi_{0,n,J_1}\circ\cdots\circ\Phi_{0,n,1}.
$
Then
\[
F_{0,n}(Y)=E_0^{J_1}Y+\Delta\widetilde W_n^{(0)},\qquad
\Delta\widetilde W_n^{(0)}
:=\sum_{j=1}^{J_1}E_0^{J_1-j}\Delta\widetilde W_{n,j}^{(0)}.
\]
The accumulated increment $\Delta\widetilde W_n^{(0)}$ belongs to $L^2(\Omega;H)$ because $J_1$ is finite and $E_0$ is bounded. In all propagator evaluations below, the same realization of the coupled Brownian increments is used across levels.

The iteration is initialized by a sequential Level-1 solve:
\[
X_{1,0}^{(0)}:=0,\qquad
X_{1,n+1}^{(0)}:=G_{1,n}(X_{1,n}^{(0)}),
\qquad 0\le n<N_1.
\]
At each outer iteration, we construct the Level-1 FAS equations from the discrepancy between the Level-0 and Level-1 propagators. These equations are approximately solved through a Level-2 FAS correction, after which F-relaxation updates the remaining Level-1 points. Let $k\in\N_0$ and suppose that $X_1^{(k)}:=(X_{1,n}^{(k)})_{n=0}^{N_1}$ is given.

For $0\le n<N_1$, define the Level-1 FAS correction
\[
\begin{aligned}
	\tau_{1,n}^{(k)}:=F_{0,n}(X_{1,n}^{(k)})-G_{1,n}(X_{1,n}^{(k)})=(E_0^{J_1}-S_1)X_{1,n}^{(k)}
	+\Delta\widetilde W_n^{(0)}-\Delta\widetilde W_n^{(1)}.
\end{aligned}
\]
The corrected Level-1 one-step propagator is
\[
\begin{aligned}
	\Psi_{1,n}^{(k)}(Y):=G_{1,n}(Y)+\tau_{1,n}^{(k)}
	=S_1Y+(E_0^{J_1}-S_1)X_{1,n}^{(k)}
	+\Delta\widetilde W_n^{(0)}.
\end{aligned}
\]
In particular, the FAS consistency identity
$
\Psi_{1,n}^{(k)}(X_{1,n}^{(k)})
=F_{0,n}(X_{1,n}^{(k)})
$
holds at the current iterate. The map $\Psi_{1,n}^{(k)}$ is an iteration-dependent random affine map with linear part $S_1$.

To compare $J_2$ corrected Level-1 steps with a single Level-2 step, define $
F_{1,m}^{(k)}
:=\Psi_{1,mJ_2+J_2-1}^{(k)}\circ\cdots\circ
\Psi_{1,mJ_2}^{(k)}, 0\le m<N_2.
$
Since each corrected Level-1 propagator has linear part $S_1$, its composition is an iteration-dependent affine map with linear part $S_1^{J_2}$. Therefore,
$
F_{1,m}^{(k)}(Y)-F_{1,m}^{(k)}(\widehat Y)
=S_1^{J_2}(Y-\widehat Y)
$
for all $Y,\widehat Y\in L^2(\Omega;H)$.

The Level-2 correction is injected into the Level-1 C-points. For sequences $V=(V_n)_{n=0}^{N_1}$ and $Z=(Z_m)_{m=0}^{N_2}$ with entries in $L^2(\Omega;H)$, define the injection restriction and prolongation by
$
(R_{\mathrm I,1}V)_m:=V_{mJ_2}, 0\le m\le N_2,
$
and $
(P_{\mathrm I,1}Z)_{mJ_2}:=Z_m,
(P_{\mathrm I,1}Z)_{mJ_2+s}:=0,
$
for $0\le m\le N_2$ in the first relation and $0\le m<N_2$, $1\le s<J_2$ in the second. Thus prolongation injects a correction only at C-points; the F-point values are subsequently updated by F-relaxation.

For $0\le m<N_2$, define the Level-2 FAS correction
$
\tau_{2,m}^{(k)}
:=F_{1,m}^{(k)}(X_{1,mJ_2}^{(k)})
-G_{2,m}(X_{1,mJ_2}^{(k)}).
$
Starting from $X_{2,0}^{(k+1)}:=0$, solve the Level-2 FAS problem sequentially $
X_{2,m+1}^{(k+1)}
:=G_{2,m}(X_{2,m}^{(k+1)})+\tau_{2,m}^{(k)},
 0\le m<N_2.$
Write $X_2^{(k+1)}:=(X_{2,m}^{(k+1)})_{m=0}^{N_2}$. The injection correction gives
\[
\begin{aligned}
	X_{1,mJ_2}^{(k+1)}
	:=X_{1,mJ_2}^{(k)}
	+\bigl(P_{\mathrm I,1}(X_2^{(k+1)}-R_{\mathrm I,1}X_1^{(k)})\bigr)_{mJ_2}
	=X_{2,m}^{(k+1)},\qquad 0\le m\le N_2.
\end{aligned}
\]
Consequently, the C-point update satisfies
\begin{equation}
	\begin{aligned}
		X_{1,(m+1)J_2}^{(k+1)}
		={}G_{2,m}(X_{1,mJ_2}^{(k+1)})
		+F_{1,m}^{(k)}(X_{1,mJ_2}^{(k)})
		-G_{2,m}(X_{1,mJ_2}^{(k)}),
	\end{aligned}
	\label{eq:level1-cpoint-recursion}
\end{equation}
for $0\le m<N_2$.

The remaining Level-1 F-point values are updated by F-relaxation. For $0\le m<N_2$ and $0\le s\le J_2-2$, set
\begin{equation}
	X_{1,mJ_2+s+1}^{(k+1)}
	=\Psi_{1,mJ_2+s}^{(k)}
	\bigl(X_{1,mJ_2+s}^{(k+1)}\bigr).
	\label{eq:level1-fpoint-recursion}
\end{equation}
Equations~\eqref{eq:level1-cpoint-recursion}--\eqref{eq:level1-fpoint-recursion} define the $(k+1)$st Level-1 iterate. The Level-0 interval propagations and the composite corrected Level-1 propagations can be evaluated concurrently across their respective intervals once the previous iterate is available. The Level-2 solve remains sequential, whereas the subsequent Level-1 F-relaxation can be carried out independently on distinct Level-2 intervals. This describes the dependency structure of the iteration; actual parallel runtimes are not assessed here.

We next define the reference trajectory against which the iteration error is measured. Set $X_{0,0,0}^{\mathrm{ref}}:=0$ and, for $0\le n<N_1$ and $1\le j\le J_1$, compute
$
X_{0,n,j}^{\mathrm{ref}}
:=\Phi_{0,n,j}(X_{0,n,j-1}^{\mathrm{ref}}).
$
At the end of each Level-1 interval, identify
$
X_{0,n+1,0}^{\mathrm{ref}}
:=X_{0,n,J_1}^{\mathrm{ref}},\qquad 0\le n<N_1,
$
and define the reference values on the Level-1 grid by
$
X_{1,n}^{\mathrm{ref}}:=X_{0,n,0}^{\mathrm{ref}}, 0\le n\le N_1.
$
It follows that $
X_{1,n+1}^{\mathrm{ref}}
=F_{0,n}(X_{1,n}^{\mathrm{ref}}),
 0\le n<N_1.
$
The reference trajectory is the Level-0 exponential Euler solution restricted to Level-1 points, not the exact mild solution of~\eqref{eq:model}.

The FAS construction is consistent with this reference trajectory. Indeed, suppose that an outer iteration has reached a fixed point $X_1^*$, so that the input and output Level-1 sequences coincide. Then the C-point update~\eqref{eq:level1-cpoint-recursion} implies
\[
X_{1,(m+1)J_2}^*=F_{1,m}^*(X_{1,mJ_2}^*),
\qquad 0\le m<N_2,
\]
where $F_{1,m}^*$ denotes the composite corrected propagator formed from $X_1^*$. Together with the F-point update~\eqref{eq:level1-fpoint-recursion}, this shows that
$X_{1,n+1}^*=\Psi_{1,n}^*(X_{1,n}^*)$ for every $0\le n<N_1$. By FAS consistency,
$X_{1,n+1}^*=F_{0,n}(X_{1,n}^*)$. Since $X_{1,0}^*=0$, uniqueness of the sequential Level-0 recursion yields $X_{1,n}^*=X_{1,n}^{\mathrm{ref}}$ for all $0\le n\le N_1$. Thus the fixed point of the three-level iteration coincides with the restricted Level-0 reference trajectory.

For $k\in\N_0$, define the Level-1 iteration error by
$
e_{1,n}^{(k)}:=X_{1,n}^{(k)}-X_{1,n}^{\mathrm{ref}},
 0\le n\le N_1.
$
The discrete noise increments are square integrable and measurable with respect to the filtration at the right endpoints of their respective temporal intervals. Using the boundedness of $E_0$, $S_1$, and $S_2$, together with the causal dependence of the FAS corrections on the preceding iterate, induction over the iteration count and temporal indices gives $
	X_{0,n,j}^{\mathrm{ref}}\in L^2(\Omega,\mathcal F_{t_{n,j}};H),$
	$X_{1,n}^{(k)},\ X_{1,n}^{\mathrm{ref}}\in L^2(\Omega,\mathcal F_{t_n};H),$ and
$X_{2,m}^{(k+1)}\in L^2(\Omega,\mathcal F_{t_{mJ_2}};H)$
for all admissible indices. In particular, $e_{1,n}^{(k)}\in L^2(\Omega,\mathcal F_{t_n};H)$, so the mean-square error quantities used in the following sections are well defined.

\section{Modal error analysis and fixed-grid convergence}
\label{sec:modal}
This section develops the error analysis of the three-level F-MGRIT
iteration relative to the Level-0 reference trajectory. We first
identify the modal propagation factors and quantify the discrepancies
between adjacent levels. We then establish the initialization error
bound and derive the coupled Hilbert-space and modal error recursions.
These identities yield a uniform modal contraction estimate while
retaining the joint dependence on the second coarsening factor and
the modal time scale. Summation over the modes gives a geometric
mean-square error bound. Finally, the causal structure of the
iteration implies finite termination on each fixed temporal grid.

\subsection[Modal representation and interlevel estimates]{Modal representation and interlevel estimates}

The deterministic operators $E_0$, $S_1$, and $S_2$ act diagonally in the common eigenbasis. 
For each $p\ge1$, set $
z_p:=\eta_p\Delta T_1,
d_p:=\mathrm e^{-\eta_p\delta t_0},
f_p:=d_p^{J_1}=\mathrm e^{-z_p},
$
and
$
b_p:=(1+z_p)^{-1},
c_p:=(1+J_2z_p)^{-1}.
$
These are the modal propagation factors associated
with the Level-0, Level-1 and Level-2 schemes.
Since $\eta_p>0$, all four factors belong to $(0,1)$.

The modal components of the reference trajectory are $x_{0,n,j,p}^{\mathrm{ref}}:=(X_{0,n,j}^{\mathrm{ref}},e_p)_H$ and $x_{1,n,p}^{\mathrm{ref}}:=(X_{1,n}^{\mathrm{ref}},e_p)_H$. For $k\in\N_0$, write $x_{1,n,p}^{(k)}:=(X_{1,n}^{(k)},e_p)_H$ for the modal component of the $k$th iterate and define $e_{1,n,p}^{(k)}:=(e_{1,n}^{(k)},e_p)_H=x_{1,n,p}^{(k)}-x_{1,n,p}^{\mathrm{ref}}$. 
For later use, set $M_p^{(k)}:=\max\limits_{0\le n\le N_1}(\E|e_{1,n,p}^{(k)}|^2)^{\frac12}$. We first derive the modal recursion satisfied by the Level-0 reference trajectory.

\begin{lemma}
	For every $p\ge1$ and \mbox{$0\le n\le N_1-1$}, the modal components of the Level-0 reference trajectory satisfy
	\begin{equation}
		x_{0,n,j,p}^{\mathrm{ref}}
		=d_p x_{0,n,j-1,p}^{\mathrm{ref}}
		+\sqrt{\gamma_p}\,d_p\Delta\beta_{n,j,p}^{(0)},
		\qquad j=1,\ldots,J_1.
		\label{eq:level0-fine}
	\end{equation}
	Composing $J_1$ Level-0 steps yields the recursion at the Level-1 points,
	\begin{equation}
		x_{1,n+1,p}^{\mathrm{ref}}
		=f_p x_{1,n,p}^{\mathrm{ref}}+\zeta_{n,p}^{(0)},
		\label{eq:level0-coarse}
	\end{equation}
	with $x_{1,0,p}^{\mathrm{ref}}=0$. Here $\zeta_{n,p}^{(0)}:=\sum_{j=1}^{J_1}\sqrt{\gamma_p}\,d_p^{J_1-j+1}\Delta\beta_{n,j,p}^{(0)}$ is the $p$th modal component of the accumulated Level-0 noise increment on $[t_n,t_{n+1}]$.
\end{lemma}

\begin{proof}
	Taking the $p$th modal component of the Level-0 reference recursion defined in Section~\ref{sec:discretization} and using
	$E_0e_p=d_pe_p$ yields~\eqref{eq:level0-fine}.
	Iterating~\eqref{eq:level0-fine}, we obtain
	\[
	x_{0,n,j,p}^{\mathrm{ref}}
	=d_p^j x_{1,n,p}^{\mathrm{ref}}
	+\sqrt{\gamma_p}\sum_{q=1}^{j}
	d_p^{j-q+1}\Delta\beta_{n,q,p}^{(0)},
	\qquad 1\le j\le J_1.
	\]
	Setting $j=J_1$ in the preceding identity and using
	$x_{1,n+1,p}^{\mathrm{ref}}=x_{0,n,J_1,p}^{\mathrm{ref}}$, we have
	\[
	x_{1,n+1,p}^{\mathrm{ref}}
		=d_p^{J_1}x_{1,n,p}^{\mathrm{ref}}
		+\sum_{q=1}^{J_1}\sqrt{\gamma_p}\,
		d_p^{J_1-q+1}\Delta\beta_{n,q,p}^{(0)}
		=f_p x_{1,n,p}^{\mathrm{ref}}+\zeta_{n,p}^{(0)}.
	\]
	The initial condition follows from $X_{1,0}^{\mathrm{ref}}=0$.
\end{proof}

The two interlevel propagation differences are $f_p-b_p$ and
$b_p^{J_2}-c_p$. To investigate robustness with respect to the
second coarsening factor, we retain the joint dependence on
$J_2$ and $z_p=\eta_p\Delta T_1$.

\begin{lemma}\label{lem:factor}
	For $z>0$ and an integer $J\ge2$, let
	$f=\mathrm{e}^{-z}$, $b=(1+z)^{-1}$, and $c=(1+Jz)^{-1}$, and define
	$\alpha(z):=\frac{b-f}{1-b}$ and $\Theta_J(z):=\frac{c-b^J}{1-c}$.
	Then $\alpha(z),\Theta_J(z)\ge0$ and
	\begin{equation}\label{eq:factor-diff}
		\begin{aligned}
			|f-b|&\le C(1\wedge z^2),
			&\alpha(z)&\le\min\left\{\frac z2,\frac1{\mathrm{e}}\right\},\\
			|b^J-c|&\le C\min\{J(J-1)z^2,1\},
			&\Theta_J(z)&\le
			\min\left\{\frac{J-1}{2}z,\frac1{Jz}\right\}.
		\end{aligned}
	\end{equation}
	The constant $C$ is independent of $z$ and $J$. Moreover,
	\begin{equation}\label{eq:theta-uniform}
		0\le\Theta_J(z)\le\alpha(Jz)\le\frac1{\mathrm{e}}.
	\end{equation}
\end{lemma}
\begin{proof}
	Since $\mathrm{e}^z\ge1+z$, we have $f\le b$. Direct calculation gives $
	\alpha(z)
	=\frac{1-(1+z)\mathrm{e}^{-z}}z
	=\frac1z\int_0^z s\mathrm{e}^{-s}\,\mathrm ds.
$
	The inequalities $s\mathrm{e}^{-s}\le s$ and $s\mathrm{e}^{-s}\le\frac1{\mathrm{e}}$ for $s\ge0$
	yield $\alpha(z)\le\min\bigl\{\frac z2,\frac1{\mathrm{e}}\bigr\}$.
	Since \mbox{$b-f=(1-b)\alpha(z)$} and $1-b\le1\wedge z$,
	the bound for $|f-b|$ follows.

	Bernoulli's inequality gives $(1+z)^J\ge1+Jz$, so $b^J\le c$.
	Using the integral remainder formula, we obtain
	\[
	\Theta_J(z)
		=\frac{(1+z)^J-1-Jz}{Jz(1+z)^J}
		=\frac{J-1}{z}\int_0^z
		(z-s)\frac{(1+s)^{J-2}}{(1+z)^J}\,\mathrm ds
		\le\frac{J-1}{2}z.
	\]
	The inequalities $0\le c-b^J\le c$ and $\frac{c}{1-c}=\frac1{Jz}$
	give $\Theta_J(z)\le\frac{1}{Jz}$.
	On the other hand, $\log(1+z)\le z$ implies
	$(1+z)^{-J}\ge \mathrm{e}^{-Jz}$, and hence
	\[
	\Theta_J(z)
	=\frac{1-(1+Jz)(1+z)^{-J}}{Jz}
	\le\frac{1-(1+Jz)\mathrm{e}^{-Jz}}{Jz}
	=\alpha(Jz)\le\frac1{\mathrm{e}}.
	\]
	Finally, $c-b^J=(1-c)\Theta_J(z)$ and $1-c\le Jz$
	yield $c-b^J\le\frac{J(J-1)}{2}z^2$. Combining this with
	$c-b^J\le1$ proves the remaining bound.
\end{proof}

\begin{remark}\label{rem:theta-regimes}
	Lemma~\ref{lem:factor} provides two complementary
	types of estimates. The bounds
	$\alpha(z),\Theta_J(z)\le\mathrm e^{-1}$
	are uniform in the modal time scale and
	the coarsening factor, and will be used
	to establish geometric convergence.
	The refined estimates $
	\alpha(z)\le \frac z2$ and
	$\Theta_J(z)\le\frac{J-1}{2}z$
	retain the dependence on the modal time scale
	and are essential for the temporal refinement
	analysis in Section~\ref{sec:refinement}.
\end{remark}
Before estimating an outer iteration, we bound the error of the initial Level-1 solve relative to the Level-0 reference trajectory. Expressing both trajectories in terms of the same Level-0 Brownian increments allows this error to be estimated by It\^o's isometry.

\begin{lemma}\label{lem:init}
	For every $p\ge1$ and $0\le n\le N_1-1$, the
	modal components of the initial Level-1 iterate
	satisfy
	\begin{equation}
		x_{1,n+1,p}^{(0)}
		=b_p x_{1,n,p}^{(0)}
		+\sqrt{\gamma_p}\,b_p\Delta\beta_{n,p}^{(1)},
		\qquad x_{1,0,p}^{(0)}=0.
		\label{eq:init}
	\end{equation}
	Moreover, there exists a constant
	$C_{\mathrm{init}}>0$, independent of
	$p,J_1,J_2,T$ and all temporal grid parameters,
	such that
	\begin{equation}
		M_p^{(0)}
		\le C_{\mathrm{init}}
		\left(\frac{\gamma_p}{\eta_p}\right)^{1/2}
		(1\wedge z_p),
		\qquad z_p=\eta_p\Delta T_1.
		\label{eq:init-est}
	\end{equation}
\end{lemma}

\begin{proof}
	The initialization in
	Section~\ref{sec:discretization} immediately
	gives~\eqref{eq:init}. Iterating this recursion
	and using the coupling of Brownian increments
	across temporal levels, we obtain
	\begin{equation}
		\begin{aligned}
			x_{1,n,p}^{(0)}
			&=\sqrt{\gamma_p}
			\sum_{i=1}^{n}\sum_{s=1}^{J_1}
			b_p^i
			\Delta\beta_{n-i,J_1-s+1,p}^{(0)},\\
			x_{1,n,p}^{\mathrm{ref}}
			&=\sqrt{\gamma_p}
			\sum_{i=1}^{n}\sum_{s=1}^{J_1}
			f_p^{i-1}d_p^s
			\Delta\beta_{n-i,J_1-s+1,p}^{(0)}.
		\end{aligned}
		\label{eq:init-representations}
	\end{equation}
	Consequently, It\^o's isometry yields
	\begin{equation}
		\mathbb E|e_{1,n,p}^{(0)}|^2
		=\gamma_p\delta t_0
		\sum_{i=1}^{n}\sum_{s=1}^{J_1}
		\left|b_p^i-f_p^{i-1}d_p^s\right|^2.
		\label{eq:init-error-isometry}
	\end{equation}
	
	Using the decomposition $	b_p^i-f_p^{i-1}d_p^s
	=(b_p^i-f_p^i)
	+f_p^{i-1}(f_p-d_p^s),$
	we obtain
	\begin{equation}
		\mathbb E|e_{1,n,p}^{(0)}|^2
		\le \mathcal I_{1,p}^{(n)}
		+\mathcal I_{2,p}^{(n)},
		\label{eq:init-error-decomposition}
	\end{equation}
	where $
		\mathcal I_{1,p}^{(n)}
		:=2\gamma_p\Delta T_1
		\sum_{i=1}^{n}|b_p^i-f_p^i|^2$ and
		$\mathcal I_{2,p}^{(n)}
		:=2\gamma_p\delta t_0
		\sum_{i=1}^{n}f_p^{2i-2}
		\sum_{s=1}^{J_1}|f_p-d_p^s|^2.$
		
	We estimate these terms separately.
	First, suppose that $0<z_p\le1$.
	By Lemma~\ref{lem:factor},
	$|b_p-f_p|\le Cz_p^2$.
	Since $
	b_p^i-f_p^i
	=(b_p-f_p)
	\sum_{j=0}^{i-1}b_p^{i-1-j}f_p^j$
	and
	$\max\{b_p,f_p\}\le \mathrm{e}^{-z_p/2}$,
	geometric-series estimates give $
	\sum_{i\ge1}|b_p^i-f_p^i|^2
	\le Cz_p.
$
	Furthermore,
	$|f_p-d_p^s|\le z_p$ and
	$\sum_{i\ge1}f_p^{2i-2}\le Cz_p^{-1}$.
	Recalling that
	$\Delta T_1=z_p/\eta_p$, we obtain
	\begin{equation}
		\mathcal I_{1,p}^{(n)}
		+\mathcal I_{2,p}^{(n)}
		\le C\frac{\gamma_p}{\eta_p}z_p^2,
		\qquad 0<z_p\le1.
		\label{eq:init-small-z}
	\end{equation}
	
	Next, suppose that $z_p\ge1$.
	Since $b_p\le(1+z_p)^{-1}$ and
	$f_p=\mathrm{e}^{-z_p}$, geometric summation gives
	$
	\sum_{i\ge1}|b_p^i-f_p^i|^2
	\le Cz_p^{-2}.$
	Moreover,
$
	\sum_{i\ge1}f_p^{2i-2}\le C,$ 
$	\delta t_0\sum_{s=1}^{J_1}d_p^{2s}
	\le\frac{1}{2\eta_p},
$
	and
	$\Delta T_1 f_p^2\le C\eta_p^{-1}$.
	Using
	$|f_p-d_p^s|^2\le2f_p^2+2d_p^{2s}$,
	we conclude that
	\begin{equation}
		\mathcal I_{1,p}^{(n)}
		+\mathcal I_{2,p}^{(n)}
		\le C\frac{\gamma_p}{\eta_p},
		\qquad z_p\ge1.
		\label{eq:init-large-z}
	\end{equation}
	
	Combining~\eqref{eq:init-error-decomposition},
	\eqref{eq:init-small-z}, and
	\eqref{eq:init-large-z}, we obtain $
	\mathbb E|e_{1,n,p}^{(0)}|^2
	\le C\frac{\gamma_p}{\eta_p}
	(1\wedge z_p^2), 0\le n\le N_1.
$
	Taking the maximum over $n$ and then
	square roots proves~\eqref{eq:init-est}.
\end{proof}

\subsection[Coupled error propagation and one-iteration estimates]{Coupled error propagation and one-iteration estimates}

To compare the three-level iteration with an exact Level-1 FAS solve, we introduce an auxiliary trajectory for the error analysis. Given the $k$th iterate, set $\widetilde X_{1,0}^{(k+1)}:=0$ and define
$\widetilde X_{1,n+1}^{(k+1)}:=\Psi_{1,n}^{(k)}\bigl(\widetilde X_{1,n}^{(k+1)}\bigr)$
for \mbox{$0\le n\le N_1-1$}.
Its error relative to the Level-0 reference trajectory is $\widetilde e_{1,n}^{(k+1)}:=\widetilde X_{1,n}^{(k+1)}-X_{1,n}^{\mathrm{ref}}$ for \mbox{$0\le n\le N_1$}. At the Level-2 grid points $t_{mJ_2}$, define the deviation from the auxiliary trajectory by \mbox{$h_{2,m}^{(k+1)}:=X_{1,mJ_2}^{(k+1)}-\widetilde X_{1,mJ_2}^{(k+1)}$} for $0\le m\le N_2$. The following recursions describe the coupling between this deviation and the auxiliary error.

\begin{lemma}\label{lem:coupled-errors}
	For every $k\in\N_0$ and \mbox{$0\le n\le N_1-1$}, the auxiliary Level-1 error satisfies
	\begin{equation}
		\widetilde e_{1,n+1}^{(k+1)}
		=S_1\widetilde e_{1,n}^{(k+1)}
		+(E_0^{J_1}-S_1)e_{1,n}^{(k)},
		\label{eq:pred-error-H}
	\end{equation}
	with $\widetilde e_{1,0}^{(k+1)}=0$. At the Level-2 points, for $0\le m\le N_2-1$,
	\begin{equation}
		h_{2,m+1}^{(k+1)}
		=S_2h_{2,m}^{(k+1)} +(S_1^{J_2}-S_2)
		\bigl(e_{1,mJ_2}^{(k)}
		-\widetilde e_{1,mJ_2}^{(k+1)}\bigr),
		\label{eq:h-error-H}
	\end{equation}
	with $h_{2,0}^{(k+1)}=0$. For $0\le m\le N_2-1$ and $0\le s\le J_2-1$, the Level-1 error decomposes as
	\begin{equation}
		e_{1,mJ_2+s}^{(k+1)}
		=\widetilde e_{1,mJ_2+s}^{(k+1)}
		+S_1^sh_{2,m}^{(k+1)}.
		\label{eq:final-error-H}
	\end{equation}
	At the final time $t_{N_1}=T$, $e_{1,N_1}^{(k+1)}=\widetilde e_{1,N_1}^{(k+1)}+h_{2,N_2}^{(k+1)}$.

	For $p\ge1$, define
	$\widetilde e_{1,n,p}^{(k+1)}:=(\widetilde e_{1,n}^{(k+1)},e_p)_H$
	for \mbox{$0\le n\le N_1$} and $h_{2,m,p}^{(k+1)}:=(h_{2,m}^{(k+1)},e_p)_H$
	for $0\le m\le N_2$.
	Taking inner products with $e_p$, we obtain, for \mbox{$0\le n\le N_1-1$},
	\begin{equation}
		\widetilde e_{1,n+1,p}^{(k+1)}
		=b_p\widetilde e_{1,n,p}^{(k+1)}
		+(f_p-b_p)e_{1,n,p}^{(k)},
		\label{eq:pred-error}
	\end{equation}
	with $\widetilde e_{1,0,p}^{(k+1)}=0$. For $0\le m\le N_2-1$,
	\begin{equation}
		h_{2,m+1,p}^{(k+1)}
		=c_ph_{2,m,p}^{(k+1)} +(b_p^{J_2}-c_p)
		\bigl(e_{1,mJ_2,p}^{(k)}
		-\widetilde e_{1,mJ_2,p}^{(k+1)}\bigr),
		\label{eq:h-error}
	\end{equation}
	with $h_{2,0,p}^{(k+1)}=0$. Moreover, for $0\le m\le N_2-1$ and $0\le s\le J_2-1$,
	\begin{equation}
		e_{1,mJ_2+s,p}^{(k+1)}
		=\widetilde e_{1,mJ_2+s,p}^{(k+1)}
		+b_p^sh_{2,m,p}^{(k+1)}.
		\label{eq:final-error}
	\end{equation}
	At the final time $t_{N_1}=T$, $e_{1,N_1,p}^{(k+1)}=\widetilde e_{1,N_1,p}^{(k+1)}+h_{2,N_2,p}^{(k+1)}$.
\end{lemma}

\begin{proof}
	By the definitions of $\Psi_{1,n}^{(k)}$ and the reference trajectory, we have $
	\widetilde X_{1,n+1}^{(k+1)}
		=S_1\widetilde X_{1,n}^{(k+1)}
		+(E_0^{J_1}-S_1)X_{1,n}^{(k)}
		+\Delta\widetilde W_n^{(0)},$ and
		$X_{1,n+1}^{\mathrm{ref}}
		=E_0^{J_1}X_{1,n}^{\mathrm{ref}}
		+\Delta\widetilde W_n^{(0)}.$
	Subtracting the reference recursion from the auxiliary recursion yields~\eqref{eq:pred-error-H}.
	For $Y,\widehat Y\in L^2(\Omega;H)$, we have
	$G_{2,m}(Y)-G_{2,m}(\widehat Y)=S_2(Y-\widehat Y)$.
	The affine structure also gives
	\mbox{$F_{1,m}^{(k)}(Y)-F_{1,m}^{(k)}(\widehat Y)=S_1^{J_2}(Y-\widehat Y)$}.

	Using \mbox{$\widetilde X_{1,(m+1)J_2}^{(k+1)} =F_{1,m}^{(k)}\bigl(\widetilde X_{1,mJ_2}^{(k+1)}\bigr)$}
	in~\eqref{eq:level1-cpoint-recursion} gives
	\[
	\begin{aligned}
		h_{2,m+1}^{(k+1)}
		&=G_{2,m}\bigl(X_{1,mJ_2}^{(k+1)}\bigr)
		-G_{2,m}\bigl(X_{1,mJ_2}^{(k)}\bigr)
		+F_{1,m}^{(k)}\bigl(X_{1,mJ_2}^{(k)}\bigr)
		-F_{1,m}^{(k)}\bigl(\widetilde X_{1,mJ_2}^{(k+1)}\bigr)\\
		&=S_2\bigl(X_{1,mJ_2}^{(k+1)}-X_{1,mJ_2}^{(k)}\bigr)
		+S_1^{J_2}\bigl(X_{1,mJ_2}^{(k)}
		-\widetilde X_{1,mJ_2}^{(k+1)}\bigr)\\
		&=S_2h_{2,m}^{(k+1)}
		+(S_1^{J_2}-S_2)
		\bigl(e_{1,mJ_2}^{(k)}
		-\widetilde e_{1,mJ_2}^{(k+1)}\bigr),
	\end{aligned}
	\]
	which proves~\eqref{eq:h-error-H}.

	Fix $0\le m\le N_2-1$. At the Level-2 grid point $t_{mJ_2}$,
	the definition of $h_{2,m}^{(k+1)}$ gives
	$X_{1,mJ_2}^{(k+1)}-\widetilde X_{1,mJ_2}^{(k+1)}=h_{2,m}^{(k+1)}$.
	For $n=mJ_2+s$ with $0\le s\le J_2-2$, subtracting the defining recursion for $\widetilde X_{1,n+1}^{(k+1)}$ from~\eqref{eq:level1-fpoint-recursion} and using the linear part $S_1$ of $\Psi_{1,n}^{(k)}$ gives
	\[
	X_{1,n+1}^{(k+1)}-\widetilde X_{1,n+1}^{(k+1)}
	=\Psi_{1,n}^{(k)}\bigl(X_{1,n}^{(k+1)}\bigr)
	-\Psi_{1,n}^{(k)}\bigl(\widetilde X_{1,n}^{(k+1)}\bigr)
	=S_1\bigl(X_{1,n}^{(k+1)}-\widetilde X_{1,n}^{(k+1)}\bigr).
	\]
	Induction on $s$ shows that
	$X_{1,mJ_2+s}^{(k+1)}-\widetilde X_{1,mJ_2+s}^{(k+1)}=S_1^sh_{2,m}^{(k+1)}$
	for $0\le s\le J_2-1$.
	Using the definitions of $e_{1,mJ_2+s}^{(k+1)}$ and $\widetilde e_{1,mJ_2+s}^{(k+1)}$ gives~\eqref{eq:final-error-H}. At the final time $t_{N_1}=T$, the identity $e_{1,N_1}^{(k+1)}=\widetilde e_{1,N_1}^{(k+1)}+h_{2,N_2}^{(k+1)}$ follows directly from the definition of $h_{2,N_2}^{(k+1)}$.
	Finally, taking the inner products of~\eqref{eq:pred-error-H},~\eqref{eq:h-error-H},~\eqref{eq:final-error-H}, and the Hilbert-space identity at $t_{N_1}=T$ with $e_p$ and using $E_0^{J_1}e_p=f_pe_p$, $S_1^se_p=b_p^se_p$, and $S_2e_p=c_pe_p$ yields~\eqref{eq:pred-error},~\eqref{eq:h-error},~\eqref{eq:final-error}, and the modal terminal relation, respectively.
\end{proof}

The coupled recursions identify two distinct
sources of iteration error. The first is
governed by the Level-0--Level-1 propagator
discrepancy $E_0^{J_1}-S_1$, while the second
arises from the inexact solution of the
intermediate-level FAS problem and is governed
by $S_1^{J_2}-S_2$.
In particular, if the intermediate-level
FAS problem is solved exactly, then
$h_{2,m}^{(k+1)}=0$ for all $m$, and the
three-level error recursion reduces to
the corresponding two-level parareal recursion.

Combining the coupled error recursions with Lemma~\ref{lem:factor}
yields the following one-iteration estimate. The estimate provides both a uniform contraction factor and a frequency-dependent refinement bound, which will be used separately in the geometric convergence analysis and the fixed-iteration temporal error estimates.

\begin{lemma}\label{lem:one-step}
	For $z>0$ and an integer $J\ge2$, let
	$\mathcal K_J(z):=\alpha(z)+(1+\alpha(z))\Theta_J(z)$.
	Then, for every $k\in\N_0$ and $p\ge1$,
	\begin{equation}\label{eq:one-step-sharp}
		M_p^{(k+1)}
		\le\mathcal K_{J_2}(z_p)M_p^{(k)}.
	\end{equation}
	In particular, the universal constant $q_\ast:=2\mathrm{e}^{-1}+\mathrm{e}^{-2}\in(0,1)$ gives
	\begin{equation}\label{eq:one-step}
		M_p^{(k+1)}
		\le q_\ast(1\wedge J_2z_p)M_p^{(k)}
		\le q_\ast M_p^{(k)}.
	\end{equation}
	The bound $q_\ast$ is independent of $k$, $p$, both coarsening
	factors, and all temporal grid parameters.
\end{lemma}

\begin{proof}
	Set $\widetilde M_p^{(k+1)}
	:=\max\limits_{0\le n\le N_1}
	(\E|\widetilde e_{1,n,p}^{(k+1)}|^2)^{\frac12}$.
	Minkowski's inequality and geometric summation
	in~\eqref{eq:pred-error} give
	\[
	\widetilde M_p^{(k+1)}
	\le\frac{|f_p-b_p|}{1-b_p}M_p^{(k)}
	=\alpha(z_p)M_p^{(k)}.
	\]
	Likewise, iterating~\eqref{eq:h-error} yields
	\[
	\max_{0\le m\le N_2}
		(\E|h_{2,m,p}^{(k+1)}|^2)^{\frac12}
		\le\frac{|b_p^{J_2}-c_p|}{1-c_p}
		\bigl(M_p^{(k)}+\widetilde M_p^{(k+1)}\bigr)\le\Theta_{J_2}(z_p)
		(1+\alpha(z_p))M_p^{(k)}.
	\]
	Applying~\eqref{eq:final-error} for $n=mJ_2+s<N_1$ and the modal identity at $n=N_1$ in Lemma~\ref{lem:coupled-errors}, we combine the bounds for $\widetilde M_p^{(k+1)}$ and $h_{2,m,p}^{(k+1)}$ with $b_p^s\le1$ to obtain~\eqref{eq:one-step-sharp}.

	By~\eqref{eq:factor-diff} and~\eqref{eq:theta-uniform},
	for $z>0$ and $J\ge2$, $
	\mathcal K_J(z)
	\le\frac1{\mathrm{e}}+\left(1+\frac1{\mathrm{e}}\right)\frac1{\mathrm{e}}
	=q_\ast<1$.
	The bounds $\alpha(z)\le\frac z2$ and $\Theta_J(z)\le\frac{J-1}{2}z$ also give
	\[
	\mathcal K_J(z)
	\le\frac z2+\left(1+\mathrm{e}^{-1}\right)\frac{J-1}{2}z
	\le\frac{1+\mathrm{e}^{-1}}2Jz
	\le q_\ast Jz,
	\]
	where $\frac{1+\mathrm{e}^{-1}}2<q_\ast$.
	Thus $\mathcal K_J(z)\le q_\ast(1\wedge Jz)$,
	which proves~\eqref{eq:one-step}.
\end{proof}

\begin{cor}\label{cor:geometric}
	Under Assumptions~\ref{ass:operators}
	and~\ref{ass:noise}, the three-level F-MGRIT
	iteration converges geometrically to the
	Level-0 reference trajectory at the Level-1
	grid points. More precisely, let $
	q_\ast:=\frac{2}{\mathrm e}
	+\frac{1}{\mathrm e^2}<1.
$
	Then, for every $k\in\N_0$,
	\begin{equation}
		\left(
		\sum_{p\ge1}(M_p^{(k)})^2
		\right)^{1/2}
		\le
		q_\ast^k
		\left(
		\sum_{p\ge1}(M_p^{(0)})^2
		\right)^{1/2}
		\le
		C_{\mathrm{init}}q_\ast^k
		\mathcal R_{0,B}.
		\label{eq:modal-geometric}
	\end{equation}
	Consequently,
	\begin{equation}
		\max_{0\le n\le N_1}
		\|X_{1,n}^{(k)}-X_{1,n}^{\mathrm{ref}}\|
		_{L^2(\Omega;H)}
		\le
		C_{\mathrm{init}}q_\ast^k
		\mathcal R_{0,B}.
		\label{eq:geometric}
	\end{equation}
	The contraction factor $q_\ast$ is independent
	of the temporal grid parameters, the
	coarsening factors $J_1,J_2$, and the
	iteration count $k$.
\end{cor}
\begin{proof}
	By Lemma~\ref{lem:one-step}, for every
	$p\ge1$ and $k\in\N_0$, $
	M_p^{(k+1)}
	\le
	\mathcal K_{J_2}(z_p)M_p^{(k)}
	\le q_\ast M_p^{(k)}.
$
	Squaring, summing over $p$, and taking
	square roots yield
	\[
	\left(
	\sum_{p\ge1}(M_p^{(k+1)})^2
	\right)^{1/2}
	\le
	q_\ast
	\left(
	\sum_{p\ge1}(M_p^{(k)})^2
	\right)^{1/2}.
	\]
	Iteration gives the first inequality
	in~\eqref{eq:modal-geometric}.
	Moreover, Lemma~\ref{lem:init} implies
	\[
	\begin{aligned}
		\sum_{p\ge1}(M_p^{(0)})^2
		\le
		C_{\mathrm{init}}^2
		\sum_{p\ge1}
		\frac{\gamma_p}{\eta_p}
		(1\wedge z_p)^2
		\le
		C_{\mathrm{init}}^2
		\mathcal R_{0,B}^2,
	\end{aligned}
	\]
	which proves the second inequality
	in~\eqref{eq:modal-geometric}.
	
	Finally, Parseval's identity and
	Tonelli's theorem give, for each
	$0\le n\le N_1$,
	\[
	\begin{aligned}
		\|e_{1,n}^{(k)}\|_{L^2(\Omega;H)}^2
		=
		\sum_{p\ge1}
		\mathbb E|e_{1,n,p}^{(k)}|^2
		\le
		\sum_{p\ge1}(M_p^{(k)})^2.
	\end{aligned}
	\]
	Taking the maximum over $n$ and using
	\eqref{eq:modal-geometric}
	establishes~\eqref{eq:geometric}.
\end{proof}

\begin{remark}\label{rem:operator-norm}
	The contraction estimate
	\eqref{eq:modal-geometric} holds in the
	modal maximal mean-square norm
	\[
	\|\boldsymbol e^{(k)}\|_{\mathrm{mm}}
	:=
	\left(
	\sum_{p\ge1}
	\max_{0\le n\le N_1}
	\mathbb E|e_{1,n,p}^{(k)}|^2
	\right)^{1/2}.
	\]
	In contrast, \eqref{eq:geometric} provides
	a geometric upper bound for the error
	measured in
	$\max_{0\le n\le N_1}
	\|e_{1,n}^{(k)}\|_{L^2(\Omega;H)}$.
	The latter estimate does not, by itself,
	assert a one-step contraction with factor
	$q_\ast$ in this maximal Hilbert-space
	mean-square norm.
\end{remark}

%\begin{remark}\label{rem:operator-norm}
%	Simultaneous diagonalization and Lemma~\ref{lem:factor} give
%	\[
%	\begin{aligned}
%		\|(E_0^{J_1}-S_1)(I-S_1)^{-1}\|_{\mathcal L(H)}
%		&=\sup_{p\ge1}\alpha(z_p)\le\frac1{\mathrm{e}},\\
%		\|(S_1^{J_2}-S_2)(I-S_2)^{-1}\|_{\mathcal L(H)}
%		&=\sup_{p\ge1}\Theta_{J_2}(z_p)\le\frac1{\mathrm{e}}.
%	\end{aligned}
%	\]
%	The bound $q_\ast$ need not be optimal.
%	The contraction estimate~\eqref{eq:one-step} concerns the modal
%	quantities $M_p^{(k)}$. Summation over the modes yields the
%	geometric mean-square bound~\eqref{eq:geometric}.
%	For an unbounded spectrum,
%	$\sup_{p\ge1}(1\wedge\eta_p\Delta T_1)=1$ for every
%	$\Delta T_1>0$. Retaining the modal time-scale dependence
%	and using noise summability are therefore necessary in the
%	present argument for error estimates under temporal refinement.
%\end{remark}

\subsection[Finite termination on fixed temporal grids]{Finite termination on fixed temporal grids}

The preceding modal bounds quantify error reduction in one outer iteration. We now turn to a complementary, purely algebraic property: causality of the coupled recursions implies finite termination on each fixed temporal grid. This argument is carried out in the Hilbert space and does not require modal estimates.

\begin{theorem}
	\label{thm:fixed-grid}
	Under Assumptions~\ref{ass:operators}
	and~\ref{ass:noise}, the three-level F-MGRIT
	iteration recovers the Level-0 reference
	trajectory in finitely many iterations.
	More precisely, for every $k\in\N_0$,
	\begin{equation}
		e_{1,n}^{(k)}=0,
		\qquad
		0\le n\le\min\{k,N_1\}.
		\label{eq:fixed-grid-exactness}
	\end{equation}
	Consequently,
	\begin{equation}
		X_{1,n}^{(k)}=X_{1,n}^{\mathrm{ref}},
		\qquad
		0\le n\le N_1,\quad k\ge N_1.
		\label{eq:fixed-grid-termination}
	\end{equation}
	Equivalently, the linear error-propagation
	operator $\mathcal E_{\mathrm{it}}$ is
	nilpotent, with $
	\mathcal E_{\mathrm{it}}^{N_1}=0.
$
\end{theorem}

\begin{proof}
	Since all propagators are affine and their
	stochastic increments are coupled across
	levels, the iteration errors satisfy a
	homogeneous linear recursion. Moreover,
	the Level-1 FAS correction at time $t_n$
	depends only on the previous iterate at
	$t_n$, while the Level-2 correction and
	subsequent F-relaxation propagate
	information forward in time.
	
	In particular, the coupled error recursions
	in Lemma~\ref{lem:coupled-errors} imply that,
	for every $1\le n\le N_1$,
	$e_{1,n}^{(k+1)}$ depends linearly only on
	$\{e_{1,j}^{(k)}:0\le j<n\}$.
	Since $e_{1,0}^{(k)}=0$ for all $k$,
	the error-propagation operator on
	the noninitial Level-1 grid points
	is represented by a strictly lower
	block-triangular matrix. Hence
	$\mathcal E_{\mathrm{it}}^{N_1}=0$.
	
	The same causal structure also yields
	\eqref{eq:fixed-grid-exactness}: the first
	$k$ Level-1 grid values are exact after
	$k$ iterations, as follows immediately
	by induction on $k$. In particular,
	all Level-1 grid values coincide with
	the reference trajectory once $k\ge N_1$,
	which proves~\eqref{eq:fixed-grid-termination}.
\end{proof}

Corollary~\ref{cor:geometric} gives a geometric error bound uniform in the temporal grids, whereas Theorem~\ref{thm:fixed-grid} guarantees exact recovery on each fixed grid. Once the Level-1 values are exact, the remaining Level-0 values can be reconstructed independently on each Level-1 interval by the propagators $\Phi_{0,n,j}$.

The refined one-iteration estimate also retains a factor $1\wedge J_2z_p$, which, combined with the initialization bound, provides the frequency-dependent decay needed for the fixed-iteration temporal error analysis in the next section.

\section{Temporal error estimates for fixed iteration counts}
\label{sec:refinement}

Section~\ref{sec:modal} establishes
geometric convergence of the iteration
and finite termination on fixed temporal
grids. We now investigate the error
under temporal refinement with a fixed
number of outer iterations. By exploiting
the frequency-dependent modal estimates,
we derive convergence bounds that
explicitly account for the second
coarsening factor $J_2$ and the spectral
regularity of the noise. Throughout
this section, errors are measured
relative to the Level-0 reference
trajectory restricted to the Level-1
grid points.

\begin{theorem}\label{thm:main}
	Under Assumption~\ref{ass:operators}, let the temporal grids and
	three-level F-MGRIT iteration be defined as in
	Section~\ref{sec:discretization}. Let $k\in\N_0$ and
	$J_1,J_2\in\N$ with $J_1,J_2\ge2$.
	If $\mathcal R_{r,B}<\infty$ for some $r\in(0,k+1]$, then
	$X_{1,n}^{\mathrm{ref}},X_{1,n}^{(\nu)}\in L^2(\Omega;H)$ for
	$0\le\nu\le k$ and \mbox{$0\le n\le N_1$}. Moreover, for every $p\ge1$, the following
	modal and Hilbert-space error estimates
	\begin{equation}\label{eq:modal-k}
		M_p^{(k)}\le C_{\mathrm{init}}q_\ast^k
		\left(\frac{\gamma_p}{\eta_p}\right)^{\frac12}
		(1\wedge z_p)(1\wedge J_2z_p)^k,
		\qquad z_p=\eta_p\Delta T_1,
	\end{equation}
	and
	\begin{equation}\label{eq:main}
		\max_{0\le n\le N_1}
		\|X_{1,n}^{(k)}-X_{1,n}^{\mathrm{ref}}\|_{L^2(\Omega;H)}
		\le C_{\mathrm{init}}q_\ast^kJ_2^{(r-1)_+}
		\mathcal R_{r,B}\Delta T_1^r
	\end{equation}
hold, respectively.
\end{theorem}
\begin{proof}
	Since $\eta_p\ge\eta_0>0$ and $r>0$,
	$\mathcal R_{0,B}^2\le\eta_0^{-2r}\mathcal R_{r,B}^2<\infty$.
	Thus Assumption~\ref{ass:noise} holds. The square integrability
	of the reference trajectory and the iterates follows from
	Section~\ref{sec:discretization}.

	Repeated application of~\eqref{eq:one-step}, followed
	by~\eqref{eq:init-est}, gives
	\[
	M_p^{(k)}
		\le q_\ast^k(1\wedge J_2z_p)^kM_p^{(0)} 
		\le C_{\mathrm{init}}q_\ast^k
		\left(\frac{\gamma_p}{\eta_p}\right)^{\frac12}
		(1\wedge z_p)(1\wedge J_2z_p)^k,
	\]
	which proves~\eqref{eq:modal-k}.
	If $0<r\le1$, then
	$(1\wedge z)(1\wedge J_2z)^k\le1\wedge z\le z^r$.
	If $1<r\le k+1$, then $0<r-1\le k$ and
	$(1\wedge z)(1\wedge J_2z)^k\le z(1\wedge J_2z)^{r-1}\le J_2^{r-1}z^r$.
	Hence, for every $z>0$ and $0<r\le k+1$,
	\begin{equation}\label{eq:coarsening-robust}
		(1\wedge z)(1\wedge J_2z)^k
		\le J_2^{(r-1)_+}z^r.
	\end{equation}
	Combining~\eqref{eq:coarsening-robust} with~\eqref{eq:modal-k} yields
	\[
	\sum_{p\ge1}(M_p^{(k)})^2
		\le C_{\mathrm{init}}^2q_\ast^{2k}J_2^{2(r-1)_+}
		\sum_{p\ge1}\frac{\gamma_p}{\eta_p}
		(\eta_p\Delta T_1)^{2r} 
		=C_{\mathrm{init}}^2q_\ast^{2k}J_2^{2(r-1)_+}
		\mathcal R_{r,B}^2\Delta T_1^{2r}.
	\]
	Finally, Parseval's identity and Tonelli's theorem give
	\[
	\max_{0\le n\le N_1}
	\|e_{1,n}^{(k)}\|_{L^2(\Omega;H)}^2
	=\max_{0\le n\le N_1}\sum_{p\ge1}\E|e_{1,n,p}^{(k)}|^2
	\le\sum_{p\ge1}(M_p^{(k)})^2.
	\]
	Taking square roots proves~\eqref{eq:main}.
\end{proof}

For fixed $J_2$, the step sizes $\Delta T_1$ and
$\Delta T_2=J_2\Delta T_1$ yield the same convergence orders.
The exponent $(r-1)_+$ in~\eqref{eq:main} makes the dependence
on the second coarsening factor explicit. In particular, the
following estimate is uniform in both coarsening factors.

\begin{cor}\label{cor:coarsening-robust}
	Let Assumption~\ref{ass:operators} hold.
	Suppose that $\mathcal R_{r,B}<\infty$
	for some $r\in(0,1]$.
	Then, for every $k\in\N_0$,
	\begin{equation}\label{eq:coarsening-uniform}
		\max_{0\le n\le N_1}
		\|X_{1,n}^{(k)}-X_{1,n}^{\mathrm{ref}}\|
		_{L^2(\Omega;H)}
		\le C_{\mathrm{init}}q_\ast^k
		\mathcal R_{r,B}\Delta T_1^r.
	\end{equation}
	The bound is uniform in the coarsening
	factors $J_1$ and $J_2$.
\end{cor}

\begin{proof}
	The assertion follows from
	Theorem~\ref{thm:main}, since
	$(r-1)_+=0$ for $0<r\le1$.
\end{proof}

%\begin{cor}\label{cor:coarsening-robust}
%	Under the assumptions of Theorem~\ref{thm:main}, if
%	$0<r\le1$ and $\mathcal R_{r,B}<\infty$, then, for every $k\in\N_0$,
%	\begin{equation}\label{eq:coarsening-uniform}
%		\max_{0\le n\le N_1}
%		\|X_{1,n}^{(k)}-X_{1,n}^{\mathrm{ref}}\|_{L^2(\Omega;H)}
%		\le C_{\mathrm{init}}q_\ast^k\mathcal R_{r,B}\Delta T_1^r.
%	\end{equation}
%\end{cor}
%
%\begin{proof}
%	For $0<r\le1$, the exponent $(r-1)_+$ in~\eqref{eq:main} is zero.
%\end{proof}
	The uniformity in~\eqref{eq:coarsening-uniform} is expressed in
	terms of the Level-1 step $\Delta T_1$. The factor
	$J_2^{(r-1)_+}$ in~\eqref{eq:main} is an upper bound and is
	not asserted to be sharp. The geometric factor $q_\ast^k$
	controls decay with the iteration count, while the admissible
	temporal order is determined by $r\le k+1$ and
	$\mathcal R_{r,B}<\infty$.

We next translate the summability condition into explicit regularity
thresholds for power-law spectra.

\begin{cor}\label{cor:power-law}
	Under Assumption~\ref{ass:operators}, let the temporal grids
	and the three-level F-MGRIT iteration be defined as in
	Section~\ref{sec:discretization}.
	Suppose that there exist $a_\eta>0$, $a_\gamma\in\R$,
	$p_0\in\N$, and constants $m_\eta,M_\eta,C_\gamma>0$
	such that, for every $p\ge p_0$,
	\begin{equation}\label{eq:power-law}
		m_\eta p^{a_\eta}\le\eta_p\le M_\eta p^{a_\eta},
		\qquad
		0\le\gamma_p\le C_\gamma p^{-a_\gamma}.
	\end{equation}
	Assume that $a_\eta+a_\gamma>1$, and define $
	\rho_B:=\frac{a_\eta+a_\gamma-1}{2a_\eta}>0.
$
	Then $\mathcal R_{r,B}<\infty$ for every
	$r\in[0,\rho_B)$.
	For any fixed $k\in\N_0$ and every
	$0<r\le k+1$ satisfying $r<\rho_B$, the
	three-level F-MGRIT iterates satisfy
	\begin{equation}\label{eq:power-law-error}
		\max_{0\le n\le N_1}
		\|X_{1,n}^{(k)}-X_{1,n}^{\mathrm{ref}}\|
		_{L^2(\Omega;H)}
		\le
		C_{\mathrm{init}}q_\ast^k
		J_2^{(r-1)_+}
		\mathcal R_{r,B}\Delta T_1^r.
	\end{equation}
	In particular, the following convergence orders are obtained:
	\begin{enumerate}[label=\textup{(\roman*)}]
		\item If $k+1<\rho_B$, choosing $r=k+1$
		in~\eqref{eq:power-law-error} yields
		the temporal convergence order $k+1$.
		
		\item If $k+1\ge\rho_B$, the estimate
		\eqref{eq:power-law-error} holds for every
		$r\in(0,\rho_B)$, yielding any temporal
		convergence order strictly below $\rho_B$.
	\end{enumerate}
	The constants $C_{\mathrm{init}}$ and $q_\ast$
	are independent of the iteration count,
	the coarsening factors, and the temporal
	grid parameters.
\end{cor}

\begin{proof}
	For every $r\in[0,\rho_B)$,
	the spectral bounds~\eqref{eq:power-law}
	imply that, for all sufficiently large $p$,
$
	\eta_p^{2r-1}\gamma_p
	\le C_r p^{a_\eta(2r-1)-a_\gamma},
$
	where $C_r>0$ is independent of $p$.
	Since $r<\rho_B$ is equivalent to
$
	a_\eta(2r-1)-a_\gamma<-1,
$
	the corresponding series converges.
	Consequently,
	\[
	\mathcal R_{r,B}^2
	=\sum_{p\ge1}\eta_p^{2r-1}\gamma_p
	<\infty.
	\]
	The estimate~\eqref{eq:power-law-error}
	then follows directly from
	Theorem~\ref{thm:main}.
	The two convergence-order assertions
	are immediate consequences of the
	admissible range $0<r\le k+1$.
\end{proof}

	The upper bound on $\gamma_p$ in~\eqref{eq:power-law} guarantees
	summability for $r<\rho_B$, but does not establish failure of
	summability for $r\ge\rho_B$. If additionally
	$\gamma_p\asymp p^{-a_\gamma}$ for all sufficiently large $p$,
	then $\mathcal R_{r,B}<\infty$ if and only if $r<\rho_B$.
	This is sharpness of the spectral summability threshold, not
	necessarily sharpness of the numerical error convergence rate.
We finally specialize the preceding result to the one-dimensional
Dirichlet Laplacian and the covariance used in~\cite{BrehierWang2020}.

\begin{cor}\label{cor:Brehier-Wang}
	Under Assumption~\ref{ass:operators},
	let $H=L^2(0,1)$, $B=0$, and
	let $A$ be the one-dimensional
	Dirichlet Laplacian with
$
	Av=v'',
	D(A)=H^2(0,1)\cap H_0^1(0,1).
$
	For $\bar\alpha>0$, suppose that
$
	e_p(x)=\sqrt2\sin(\pi px),
	\lambda_p=(\pi p)^2,
	\gamma_p=\lambda_p^{1/2-2\bar\alpha}, p\ge1.
$
	Then $\mathcal R_{r,B}<\infty$
	if and only if $r<\bar\alpha$
	for $r\ge0$.
	For every $k\in\N_0$ and
	$0<r\le k+1$ satisfying
	$r<\bar\alpha$, we have
	\begin{equation}\label{eq:dirichlet-error}
		\max_{0\le n\le N_1}
		\|X_{1,n}^{(k)}-X_{1,n}^{\mathrm{ref}}\|
		_{L^2(\Omega;H)}
		\le C_{\mathrm{init}}q_\ast^k
		J_2^{(r-1)_+}
		\mathcal R_{r,B}\Delta T_1^r.
	\end{equation}
	Thus the order $k+1$ is obtained
	when $k+1<\bar\alpha$, while any
	order below $\bar\alpha$ is
	admissible when $k+1\ge\bar\alpha$.
	For $0<r\le1$, the bound is
	uniform in both coarsening factors.
\end{cor}

\begin{proof}
	Since $B=0$, $
	\eta_p=\lambda_p=(\pi p)^2,
	\gamma_p=(\pi p)^{1-4\bar\alpha}.
$
	Consequently,
	\[
	\mathcal R_{r,B}^2
	=
	\pi^{4r-1-4\bar\alpha}
	\sum_{p\ge1}p^{4r-1-4\bar\alpha},
	\]
	which is finite if and only if
	$r<\bar\alpha$.
	The power-law exponents are
	$a_\eta=2$ and
	$a_\gamma=4\bar\alpha-1$,
	so that $\rho_B=\bar\alpha$.
	The asserted estimates follow from
	Corollary~\ref{cor:power-law}.
\end{proof}

At the level of admissible temporal convergence orders,
Corollary~\ref{cor:Brehier-Wang} recovers the corresponding
regularity-dependent order in Theorem~4.1 of~\cite{BrehierWang2020}. The present estimates additionally
make explicit the dependence on the second coarsening factor
for the three-level iteration.

\section{Numerical experiments}
\label{sec:numerical}
The experiments examine the temporal error estimates and coarsening-factor dependence in Theorem~\ref{thm:main} and Corollary~\ref{cor:power-law}, together with iteration errors on fixed temporal grids. We also compare iteration counts and temporal critical-path work under a unit-cost model.

We consider $H=L^2(0,1)$ with the orthonormal basis $
e_p(x)=\sqrt{2}\sin(\pi px), p\ge1.
$
The operators $A$, $B$, and $Q$ are diagonal in this basis,
with $Ae_p=-(\pi p)^2e_p,
Be_p=e_p,$ and
$Qe_p=\gamma_pe_p.$
Consequently,
$
\eta_p=(\pi p)^2-1,
(\pi^2-1)p^2\le\eta_p\le\pi^2p^2,
 p\ge1.
$
For the numerical implementation, we project the stochastic
evolution equation onto the finite-dimensional subspace
$
H_P:=\operatorname{span}\{e_1,\ldots,e_P\},
$
and retain only the first $P$ modes of the stochastic forcing.
The initial condition is projected onto $H_P$ accordingly.
Since the drift operators are diagonal in the chosen basis,
both the Level-0 reference trajectory and the three-level
F-MGRIT iterates remain in $H_P$.
Thus, all numerical errors reported below are evaluated
for the same $P$-mode spectral approximation.

Unless otherwise specified, the temporal refinement
experiments use
$
T=1, J_1=J_2=4,
P=128$ and  $N_{\mathrm{MC}}=200
$
with $k=0,\ldots,4$ and $N_1=2^j$ for $j=4,\ldots,10$.
Accordingly, $
\Delta T_1=2^{-j},
\delta t_0=2^{-(j+2)},$ 
$\Delta T_2=4\Delta T_1$ and 
$N_2=N_1/4.
$
The temporal propagators and stochastic increments are
defined as in Section~\ref{sec:discretization}.
On each temporal grid, the Level-0 reference trajectory
and the corresponding F-MGRIT iterates are driven by
the same Brownian increments, consistently with the
interlevel coupling used in the theoretical analysis.
Within each noise-spectrum experiment, common random
numbers are used across the parameter cases.

For the $\ell$th of the $N_{\mathrm{MC}}$ Monte Carlo samples, let $X_{1,n,P}^{(k,\ell)}$ and
$X_{1,n,P}^{\mathrm{ref},(\ell)}$ denote the truncations to the first $P$ modes of the Level-1 iterate and the Level-0 reference trajectory at $t_n$, respectively.
Define
\[
  \widehat E_k(\Delta T_1)
  :=
  \max_{0\le n\le N_1}
  \left[
    \frac{1}{N_{\mathrm{MC}}}
    \sum_{\ell=1}^{N_{\mathrm{MC}}}
    \left\|
      X_{1,n,P}^{(k,\ell)}
      -X_{1,n,P}^{\mathrm{ref},(\ell)}
    \right\|_H^2
  \right]^{\frac12}.
\]

The quantity $\widehat E_k(\Delta T_1)$
estimates the maximum over Level-1 grid
points of the root-mean-square (RMS) iteration
error for the $P$-mode approximation.
The maximum is taken after averaging
the squared errors over the Monte Carlo
samples, consistently with the error
norm used in Theorem~\ref{thm:main}.

Fitted orders $r_{\mathrm{fit}}$ are obtained by least-squares regression of $\log\widehat E_k$ against $\log\Delta T_1$ on the three finest grids, $\Delta T_1\in\{2^{-8},2^{-9},2^{-10}\}$. 
A resolution check was performed
over the tested temporal grids and
iteration counts. Increasing the
number of retained modes from $96$
to $128$ changed the computed
iteration errors by at most
$1.4\times10^{-14}$ in relative terms,
indicating that the additional
high-frequency contributions were
numerically negligible for these
parameter settings.

%Preliminary resolution checks showed that increasing $N_{\mathrm{MC}}$ from $150$ to $200$ changed the finest-grid error estimates by less than $5\%$ and the fitted orders by less than $0.05$, while increasing $P$ from $96$ to $128$ produced a maximum relative change below $1.4\times10^{-14}$ in the error estimates.

\subsection{Temporal convergence}

\begin{figure}[!htbp]
  \centering
  \includegraphics[width=0.75\textwidth]{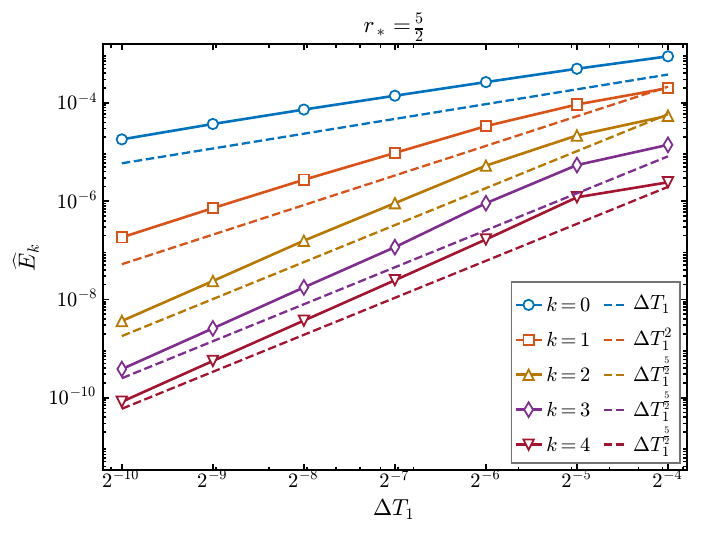}
  \caption{Temporal convergence for the logarithmically modified spectrum with $r_\ast=\frac52$.}
  \label{fig:general-convergence}
\end{figure}

We first consider a logarithmically
modified noise spectrum for which
the endpoint regularity is attained.
Let $r_\ast=5/2$ and define
\[
\gamma_p=
\eta_p^{1-2r_\ast}
\frac{1+\frac14\sin^2p}
{p\log^2(p+1)},\qquad p\ge1.
\]
Since $\eta_p\asymp p^2$, we obtain $\gamma_p\le Cp^{-9}\log^{-2}(p+1)$, so $Q$ is trace class. Moreover,
$
\eta_p^{2r-1}\gamma_p
\asymp
\frac{p^{4(r-r_\ast)-1}}
{\log^2(p+1)}.
$
Consequently,
$\mathcal R_{r,B}<\infty$
if and only if $r\le r_\ast$.
In particular, the endpoint
$r=r_\ast$ is admissible.
Set \mbox{$r_k:=\min\{k+1,r_\ast\}$}. Theorem~\ref{thm:main} provides the following estimate
\[
  \max_{0\le n\le N_1}
  \|X_{1,n}^{(k)}-X_{1,n}^{\mathrm{ref}}\|_{L^2(\Omega;H)}
  \le C_{\mathrm{init}}q_\ast^kJ_2^{(r_k-1)_+}\mathcal R_{r_k,B}\Delta T_1^{r_k}.
\]
 Fig.~\ref{fig:general-convergence} shows the corresponding temporal-refinement errors for this logarithmically modified spectrum.

The fitted orders for $k=0$ and $k=1$
are close to the predicted orders
$1$ and $2$ in Theorem~\ref{thm:main}, respectively. For
$k=2,3,4$, the observed slopes
exceed the guaranteed order $5/2$.
This is consistent with the theoretical
upper bounds, which do not assert
sharpness.
Moreover, the experiments employ
a fixed spectral truncation, for
which all finite-order spectral
regularity quantities are finite.
Consequently, the observed slopes
on the tested grids should not be
interpreted as evidence that the
infinite-dimensional regularity
threshold is sharp or can be improved.

We next examine the dependence on noise regularity in Corollary~\ref{cor:power-law} using power-law spectra. Let $\gamma_p=\eta_p^{\frac12-2\rho_B}$ for
\mbox{$\rho_B\in\{\frac32,\frac52,4,6\}$}. Since $\eta_p\asymp p^2$, the power-law parameters are $a_\eta=2$ and $a_\gamma=4\rho_B-1$, and
$\eta_p^{2r-1}\gamma_p\asymp p^{4(r-\rho_B)-1}$; hence
$\mathcal R_{r,B}<\infty$ exactly when $r<\rho_B$.
We first examine the effect of the iteration count. Fix $\rho_B\in\{\frac52,6\}$. Since the power-law spectrum does not satisfy
$\mathcal R_{r,B}<\infty$ at $r=\rho_B$, we plot a reference line of slope $2.45$, which is admissible and close to $\rho_B=\frac52$.
For the cases in which the spectral
regularity threshold limits the
admissible order, we select reference
slopes slightly below that threshold.
These slopes represent convergence
orders guaranteed by the theory,
rather than asserted sharp
asymptotic rates.

\begin{figure}[!htbp]
  \centering
  \includegraphics[width=\textwidth]{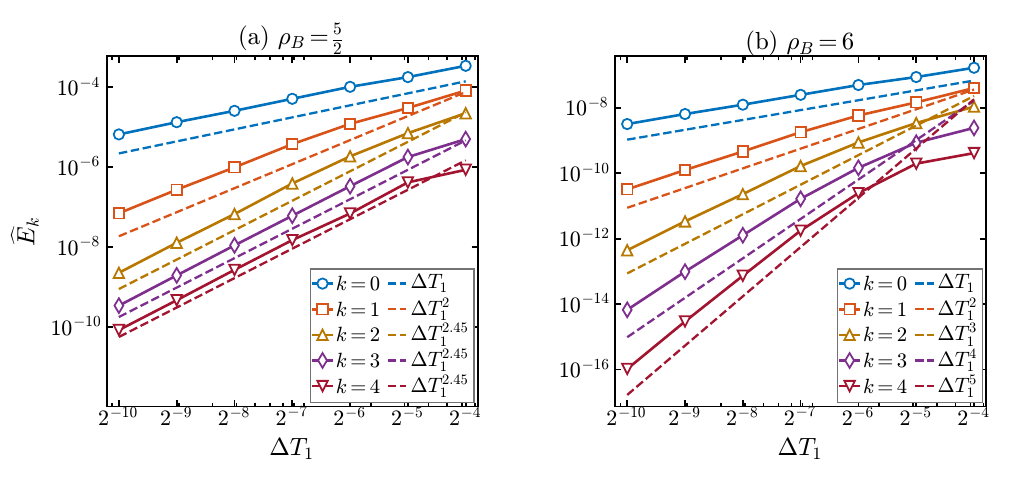}
  \caption{Temporal convergence for power-law spectra: \textup{(a)} $\rho_B=\frac52$; \textup{(b)} $\rho_B=6$.}
  \label{fig:power-vary-k}
\end{figure}
On the three finest grids, the fitted orders in Fig.~\ref{fig:power-vary-k} for $\rho_B=\frac52$ are close to the reference slopes $1$ and $2$ for $k=0,1$ and to $2.45$ for $k=2,3,4$. For $\rho_B=6$, the fitted orders for $k=0,1,2,3$ are close to the reference slopes $1$, $2$, $3$, and $4$, respectively. For $\rho_B=6$ and $k=4$,
the fitted order is $4.728$,
compared with the guaranteed order
$5$ in Corollary~\ref{cor:power-law}. The discrepancy may reflect
finite-resolution effects, and the
available grid range does not
establish the asymptotic rate.
 For each fixed $\Delta T_1$, $\widehat E_k(\Delta T_1)$ decreases over the computed iterations.

To examine the effect of the spectral regularity index, fix $k=2$ and take \mbox{$\rho_B\in\{\frac32,\frac52,4,6\}$}. Since $\mathcal R_{r,B}<\infty$ only for $r<\rho_B$, we use reference slopes $1.45$ and $2.45$ for $\rho_B=\frac32$ and $\rho_B=\frac52$, respectively; both are admissible and close to the corresponding values of $\rho_B$. For $\rho_B=4$ and $\rho_B=6$, we use slope $3=k+1$.

\begin{figure}[!ht]
  \centering
  \includegraphics[width=0.75\textwidth]{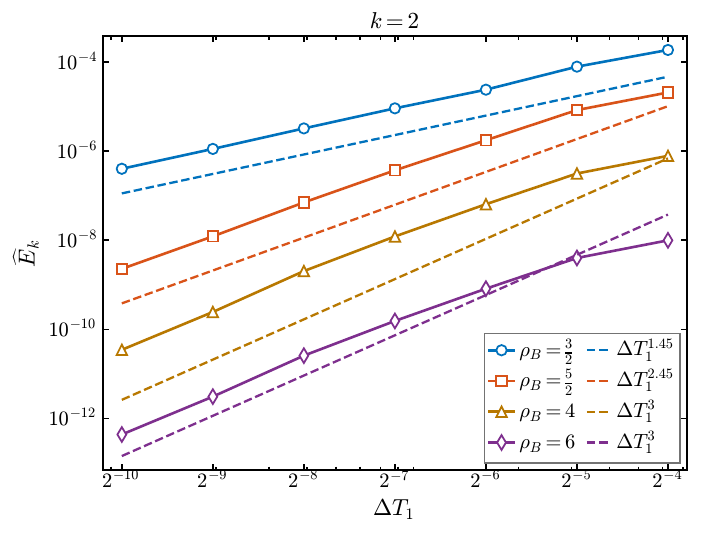}
  \caption{Temporal convergence for power-law spectra with $k=2$ and $\rho_B\in\{\frac32,\frac52,4,6\}$.}
  \label{fig:power-vary-rho}
\end{figure}

On the three finest grids, the fitted orders in Fig.~\ref{fig:power-vary-rho} are close to the reference slopes $1.45$ and $2.45$ for $\rho_B=\frac32$ and $\rho_B=\frac52$, respectively, and to $3$ for $\rho_B=4$ and $\rho_B=6$. Corollary~\ref{cor:power-law} guarantees any order $r<\rho_B$ for $\rho_B=\frac32,\frac52$ and the order $r=3=k+1$ for $\rho_B=4,6$.
Table~\ref{tab:fitted-order-summary} reports all five fitted orders for the logarithmically modified spectrum and selected fitted orders for the power-law cases $\rho_B=\frac52$ and $\rho_B=6$. The coefficient of determination exceeds $0.999$ in every row of Table~\ref{tab:fitted-order-summary}.

\begin{table}[!htbp]
  \centering
  \caption{Reference slopes and fitted orders for logarithmically modified and power-law spectra}
  \label{tab:fitted-order-summary}
  \small
  \renewcommand{\arraystretch}{1.25}
  \begin{tabular*}{\textwidth}{@{\extracolsep{\fill}}lcccc@{}}
    \toprule
    Spectrum & Parameter & Iteration $k$ & Reference slope & Fitted order $r_{\mathrm{fit}}$\\
    \midrule
    Non-power-law & $r_\ast=\frac52$ & $0$ & $1$ & $1.008$\\
    Non-power-law & $r_\ast=\frac52$ & $1$ & $2$ & $1.947$\\
    Non-power-law & $r_\ast=\frac52$ & $2$ & $2.5$ & $2.714$\\
    Non-power-law & $r_\ast=\frac52$ & $3$ & $2.5$ & $2.761$\\
    Non-power-law & $r_\ast=\frac52$ & $4$ & $2.5$ & $2.743$\\
    Power-law & $\rho_B=\frac52$ & $2$ & $2.45$ & $2.443$\\
    Power-law & $\rho_B=6$ & $3$ & $4$ & $3.782$\\
    \bottomrule
  \end{tabular*}
\end{table}

\FloatBarrier
\begingroup
\subsection{Dependence on the coarsening factors}
\label{subsec:coarsening-factors}

We examine the coarsening-factor dependence in Theorem~\ref{thm:main}
with $T=1$, \mbox{$P=128$},
$\rho_B=\frac52$, $\gamma_p=\eta_p^{-\frac92}$ and $k=2$.
For each $N_1\in\{64,128,256\}$, we vary
\mbox{$J_1\in\{2,4,8,16,32\}$} at fixed $J_2=4$, and then vary $J_2$
over the same set at fixed $J_1=4$.
The step size $\Delta T_1=\frac{1}{N_1}$ is fixed along each curve.
As $J_1$ varies, $\delta t_0=\frac{\Delta T_1}{J_1}$ and the Level-0 reference
trajectory change; each error is measured relative to the
corresponding reference trajectory.

Let $X_{1,n,P}^{(k)}$ and $X_{1,n,P}^{\mathrm{ref}}$ denote the
$P$-mode truncations of the iterate and the reference trajectory,
respectively, and set
$$
  E_{k,P}:=\max_{0\le n\le N_1}
  \|X_{1,n,P}^{(k)}-X_{1,n,P}^{\mathrm{ref}}\|_{L^2(\Omega;H)},
  \qquad
  \mathcal R_{2,B,P}:=\left(\sum_{p=1}^{P}\eta_p^3\gamma_p\right)^{\frac12}.
$$
\begin{figure}[!htbp]
  \centering
  \includegraphics[width=\textwidth]{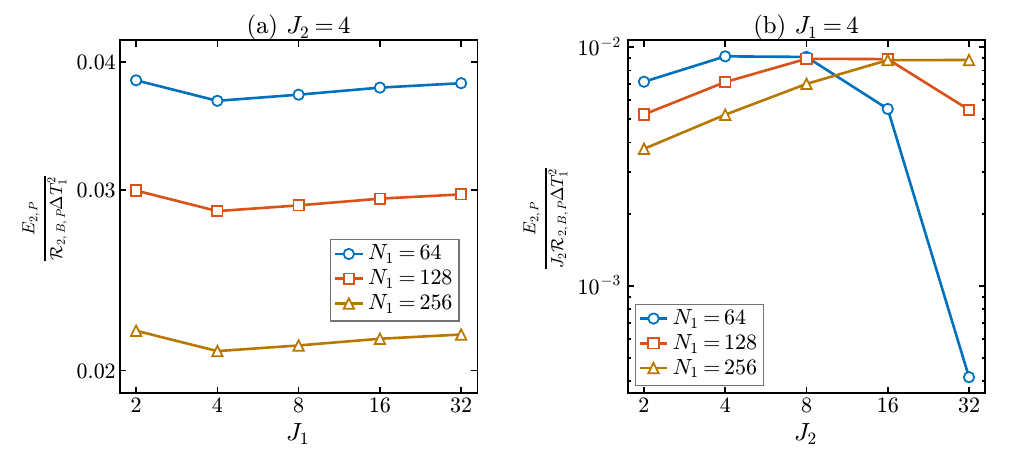}
    \caption{Normalized iteration errors for $k=2$: \textup{(a)} varying $J_1$ with $J_2=4$; \textup{(b)} varying $J_2$ with $J_1=4$.}
  \label{fig:coarsening-factors}
\end{figure}

The errors $E_{k,P}$ are evaluated
from the exact second moments of
the modal error recursions, with
the correlations induced by the
shared Brownian increments retained.
Thus, unlike the Monte Carlo
estimates used in the temporal
refinement experiments, these
quantities are free from sampling
error.
For $k=r=2$, Theorem~\ref{thm:main}
gives
$
E_{2,P}\le
C_{\mathrm{init}}q_\ast^2
J_2\mathcal R_{2,B,P}\Delta T_1^2.
$
Accordingly, we examine the
normalized errors in
Fig.~\ref{fig:coarsening-factors}
to assess their dependence on
the two coarsening factors.

%The errors $E_{k,P}$ are computed directly from the second moments
%of the modal recursions, retaining the covariance between the fine
%and coarse noise increments; no Monte Carlo approximation is used.
%Applying Theorem~\ref{thm:main} to the covariance truncated to the
%first $P$ modes with $k=r=2$ gives
%$E_{2,P}\le C_{\mathrm{init}}q_\ast^2J_2\mathcal R_{2,B,P}\Delta T_1^2$.
%The prefactor $C_{\mathrm{init}}q_\ast^2$ is independent of
%$J_1$, $J_2$, and $P$.
%Fig.~\ref{fig:coarsening-factors} shows the normalized errors
%for the two parameter studies.

For each fixed $N_1$ and $J_2=4$, the ratio of the largest to the
smallest error over the tested $J_1$ values is less than $1.048$.
For $J_1=4$, the normalized errors in
Fig.~\ref{fig:coarsening-factors}\textup{(b)} remain below
$9.2\times10^{-3}$ over all tested values of $N_1$ and $J_2$.
These observations are consistent with the coarsening-factor bound
in Theorem~\ref{thm:main}; the data do not establish sharpness
of the factor $J_2$.

\endgroup
\FloatBarrier

\subsection{Fixed-grid iteration convergence and critical-path work}
\label{subsec:iteration-work}

\begin{figure}[!htbp]
  \centering
  \includegraphics[width=0.80\textwidth]{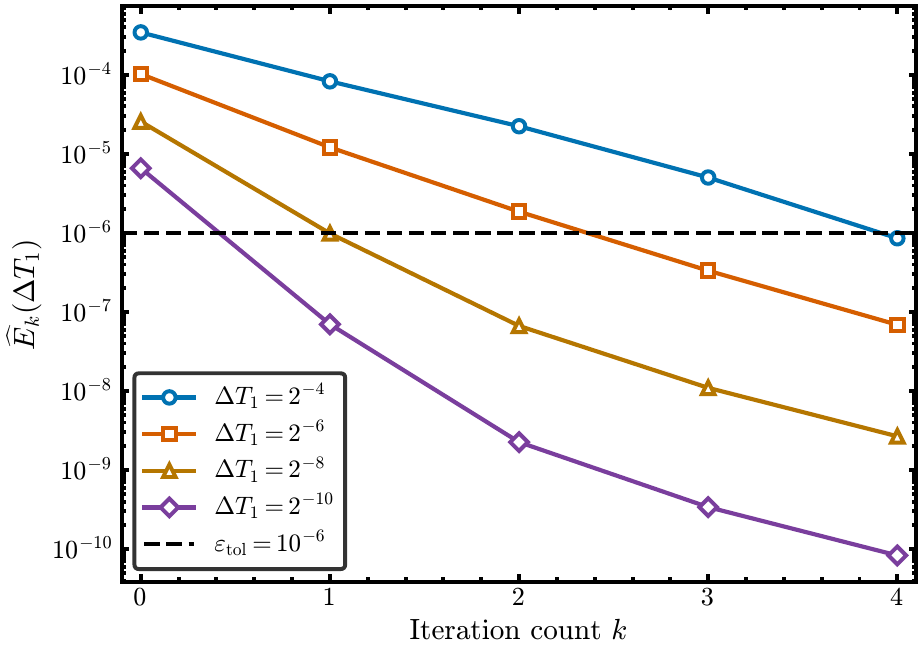}
  \caption{Fixed-grid iteration convergence for the power-law spectrum with $\rho_B=\frac52$. The dashed line marks $\varepsilon_{\mathrm{tol}}=10^{-6}$.}
  \label{fig:iteration-convergence}
\end{figure}

We next study iteration convergence on fixed temporal grids. Fix $\rho_B=\frac52$ and use
$N_1\in\{16,64,256,1024\}$, corresponding to
$\Delta T_1\in\{2^{-4},2^{-6},2^{-8},2^{-10}\}$.
With $\varepsilon_{\mathrm{tol}}=10^{-6}$, define the first computed iteration satisfying the tolerance by
\[
k_{\mathrm{tol}}(\Delta T_1)
:= \min\bigl\{k\in\{0,1,2,3,4\}: \widehat E_k(\Delta T_1)\le\varepsilon_{\mathrm{tol}}\bigr\}.
\]
Fig.~\ref{fig:iteration-convergence} shows the iteration errors on the four fixed Level-1 grids, and Table~\ref{tab:iteration-tolerance} reports the corresponding values of $k_{\mathrm{tol}}$.

\begin{table}[!htbp]
  \centering
  \caption{Iteration counts required to satisfy $\widehat E_k(\Delta T_1)\le10^{-6}$}
  \label{tab:iteration-tolerance}
  \begin{tabular*}{0.85\textwidth}{@{\extracolsep{\fill}}cccc@{}}
    \toprule
    $N_1$ & $\Delta T_1$ & $k_{\mathrm{tol}}$ & $\widehat E_{k_{\mathrm{tol}}}(\Delta T_1)$\\
    \midrule
    $16$   & $2^{-4}$  & $4$ & $8.59\times10^{-7}$\\
    $64$   & $2^{-6}$  & $3$ & $3.34\times10^{-7}$\\
    $256$  & $2^{-8}$  & $1$ & $9.89\times10^{-7}$\\
    $1024$ & $2^{-10}$ & $1$ & $7.00\times10^{-8}$\\
    \bottomrule
  \end{tabular*}
\end{table}
The errors in Fig.~\ref{fig:iteration-convergence} decrease monotonically over the five computed iterates. For the fixed absolute tolerance used here, the iteration counts in Table~\ref{tab:iteration-tolerance} do not increase under temporal refinement. These observations are consistent with the geometric mean-square bound in Corollary~\ref{cor:geometric} and the temporal error estimates in Theorem~\ref{thm:main}. The plotted Monte Carlo estimates concern the maximum-in-time error and do not directly measure the modal contraction factor in~\eqref{eq:one-step}. On each fixed grid, Theorem~\ref{thm:fixed-grid} guarantees exact recovery at all Level-1 points after at most $N_1$ iterations.

For the critical-path comparison, we use the stricter tolerance $10^{-8}$ on the fixed grid $N_1=256$. The test uses \mbox{$P=128$}, $J_1=J_2=4$, and $\rho_B=\frac52$ and compares the sequential Level-0 reference computation with parareal, two-level MGRIT (F), and three-level F-MGRIT.

The maximum-in-time root-mean-square errors are computed from the joint Gaussian covariance using~\eqref{eq:init-representations} and~\eqref{eq:pred-error}--\eqref{eq:final-error}, including the modal identity at $n=N_1$ in Lemma~\ref{lem:coupled-errors}. No Monte Carlo approximation is used; the first iteration with error at most $10^{-8}$ is reported.

For these propagators and injection operators, parareal and two-level MGRIT coincide at the Level-1 points~\cite{DobrevEtAl2017}. We assign unit cost to each single-step propagation and assume that independent interval solves are executed concurrently. With $N_0=J_1N_1$ and $N_2=\frac{N_1}{J_2}$, set $\mathcal W_{\mathrm{seq}}:=N_0$ and define
\[
\mathcal W_2(k):=N_1+k(J_1+N_1),\qquad
  \mathcal W_3(k):=N_1+k\bigl(J_1+J_2+N_2+(J_2-1)\bigr).
\]
Here $N_1$ accounts for Level-1 initialization. In $\mathcal W_3(k)$, the costs $J_1,J_2,N_2,J_2-1$ represent, respectively, the Level-0 interval solve, composite Level-1 FAS propagation, Level-2 solve, and Level-1 F-relaxation. Table~\ref{tab:efficiency-comparison} gives the error and work ratio $\frac{\mathcal W_{\mathrm{seq}}}{\mathcal W}$ for each method, where $\mathcal W$ denotes its critical-path work.

\begin{table}[!htbp]
  \centering
  \footnotesize
  \setlength{\tabcolsep}{3pt}
  \caption{Fixed-grid errors and critical-path work for tolerance $10^{-8}$}
  \label{tab:efficiency-comparison}
  \begin{tabular}{@{}lcccc@{}}
    \toprule
    Method & \shortstack{Iteration\\count} & \shortstack{Maximum-in-time\\RMS iteration\\error}
           & \shortstack{Critical-path\\work} & Work ratio\\
    \midrule
    Sequential Level-0 reference & -- & $0$ & $1024$ & $1.00$\\
    Parareal & $2$ & $7.5963\times10^{-9}$ & $776$ & $1.32$\\
    Two-level MGRIT (F) & $2$ & $7.5963\times10^{-9}$ & $776$ & $1.32$\\
    Three-level F-MGRIT & $4$ & $2.5044\times10^{-9}$ & $556$ & $1.84$\\
    \bottomrule
  \end{tabular}
\end{table}
\FloatBarrier

For the prescribed tolerance,
parareal and two-level MGRIT
reach the target accuracy in two
iterations, whereas the three-level
F-MGRIT method requires four
iterations. Under the idealized
unit-cost model, the corresponding
temporal critical-path work estimates
are $776$ and $556$, respectively.
These results illustrate a trade-off
between iteration count and
critical-path work for the tested
configurations. In particular,
the additional coarse level can
reduce the modeled critical-path
cost despite an increased number
of outer iterations. The comparison
is specific to the chosen parameters
and cost assumptions and does not
imply a general performance
advantage over parareal or
two-level MGRIT.

\FloatBarrier
\section{Conclusion}
\label{sec:conclusion}

Under the common-eigenbasis assumption, the coupled three-level error recursions yield a uniform modal contraction bound $q_\ast=2\mathrm{e}^{-1}+\mathrm{e}^{-2}<1$. Combined with the basic noise summability condition, this gives a geometric $L^2(\Omega;H)$ error bound relative to the Level-0 reference trajectory, uniformly in the temporal grids and coarsening factors. Retaining the modal time-scale dependence further gives an error bound proportional to $q_\ast^kJ_2^{(r-1)_+}\mathcal R_{r,B}\Delta T_1^r$ whenever $0<r\le k+1$ and $\mathcal R_{r,B}<\infty$, with a prefactor independent of $k$ and all temporal grid parameters. For power-law spectra, the admissible orders are $k+1$ when $k+1<\rho_B$ and any $0<r<\rho_B$ otherwise. In the specified one-dimensional Dirichlet Laplacian setting with $B=0$, these orders agree with the corresponding two-level bounds. On each fixed temporal grid, the iteration also recovers the Level-0 reference values at all Level-1 points after at most $N_1$ iterations.

The numerical experiments illustrate
the temporal error estimates and
fixed-grid iteration convergence.
A comparison with parareal and
two-level MGRIT under an idealized
unit-cost model further highlights
the trade-off between iteration
counts and estimated temporal
critical-path work in the tested
configurations.
The common-eigenbasis assumption reduces the Level-0--Level-1 and Level-1--Level-2 differences to scalar modal factors. Extending the analysis to sectorial generators without simultaneous diagonalization requires operator-level estimates in place of the present modal argument. Multiplicative noise introduces state-dependent stochastic terms into the error recursion. Assessing wall-clock performance requires parallel implementations that account for communication and synchronization costs.

\backmatter
% Before submission, add the author-confirmed Statements and Declarations
% (competing interests, data/code availability as applicable) here.
% Author contributions and competing interests must also be entered in the submission system.

\end{document}